\documentclass[11pt,twoside]{article}
\newcommand{\ifims}[2]{#1} 
\newcommand{\ifAMS}[2]{#1}   
\newcommand{\ifau}[3]{#1}  

\newcommand{\ifbook}[2]{#1}   

\ifbook{
    \newcommand{\Chapter}[1]{\section{#1}}
    \newcommand{\Section}[1]{\subsection{#1}}
    \newcommand{\Subsection}[1]{\subsubsection{#1}}

  }
  {
    \newcommand{\Chapter}[1]{\chapter{#1}}
    \newcommand{\Section}[1]{\section{#1}}
    \newcommand{\Subsection}[1]{\subsection{#1}}

  }

\def\thetitle{Bernstein - von Mises Theorem for growing parameter dimension}
\def\thanksa
{The author is partially supported by
Laboratory for Structural Methods of Data Analysis in Predictive Modeling, MIPT, 
RF government grant, ag. 11.G34.31.0073.
Financial support by the German Research Foundation (DFG) through the Collaborative 
Research Center 649 ``Economic Risk'' is gratefully acknowledged.
The author thanks the AE and two anonymous referees for very useful comments and suggestions.} 
\def\theruntitle{BvM theorem for growing parameter dimension}

\def\theabstract{
This paper revisits the prominent Fisher, Wilks, and Bernstein -- von Mises (BvM) results 
from different viewpoints.
Particular issues to address are: nonasymptotic framework with just one finite sample, 
possible model misspecification, and a large parameter dimension.
In particular,  in the case of an i.i.d. sample, the mentioned results can be stated for any smooth parametric family provided 
that the dimension \( \dimp \) of the parameter space satisfies the condition 
``\( \dimp^{2}/\nsize \) is small'' for the Fisher expansion, while the Wilks and 
the BvM results require ``\( \dimp^{3}/\nsize \) is small''.
}

\def\kwdp{62F15}
\def\kwds{62F25}

\def\thekeywords{posterior, concentration, Gaussian approximation}

\def\authora{Vladimir Spokoiny}
\def\runauthora{spokoiny, v.}
\def\addressa{
    Weierstrass-Institute, \\ Humboldt University Berlin, \\ Moscow Institute of
    Physics and Technology
    \\
    Mohrenstr. 39, 10117 Berlin, Germany,    \\
    }
\def\emaila{spokoiny@wias-berlin.de}
\def\affiliationa{Weierstrass-Institute and Humboldt University Berlin}

\renewcommand{\(}{$\,}
\renewcommand{\)}{\,$}

\def\nquad{\hspace{-1cm}}
\def\eqdef{\stackrel{\operatorname{def}}{=}}

\newcommand{\cc}[1]{\mathscr{#1}}
\newcommand{\bb}[1]{\boldsymbol{#1}}

\renewcommand{\bar}[1]{\overline{#1}}
\renewcommand{\hat}[1]{\widehat{#1}}
\renewcommand{\tilde}[1]{\widetilde{#1}}

\renewcommand{\Gamma}{\varGamma}
\renewcommand{\Pi}{\varPi}
\renewcommand{\Sigma}{\varSigma}
\renewcommand{\Delta}{\varDelta}
\renewcommand{\Lambda}{\varLambda}
\renewcommand{\Psi}{\varPsi}
\renewcommand{\Phi}{\varPhi}
\renewcommand{\Theta}{\varTheta}
\renewcommand{\Omega}{\varOmega}
\renewcommand{\Xi}{\varXi}
\renewcommand{\Upsilon}{\varUpsilon}

\def\Cov{\operatorname{Cov}}
\def\Var{\operatorname{Var}}

\def\argmax{\operatornamewithlimits{argmax}}

\def\tr{\operatorname{tr}}

\def\R{I\!\!R}
\def\E{I\!\!E}
\def\P{I\!\!P}

\def\kappa{\varkappa}

\def\T{\top}

\def\av{\bb{a}}

\def\fv{\bb{f}}

\def\uv{\bb{u}}

\def\Yv{\bb{Y}}

\def\alphav{\bb{\alpha}}
\def\epsv{\bb{\varepsilon}}
\def\etav{\bb{\eta}}
\def\gammav{\bb{\gamma}}
\def\varepsilonv{\bb{\varepsilon}}

\def\xiv{\bb{\xi}}

\def\CONST{\mathtt{C}}

\def\ND{\cc{N}}

\def\cond{\, \big| \,}

\def\rdl{\epsilon}
\def\rd{\bb{\rdl}}
\def\rddelta{\delta}
\def\rdomega{\varrho}

\def\rdb{\rd}
\def\rdm{\underline{\rdb}}

\def\Id{I\!\!\!I}
\def\Ind{\operatorname{1}\hspace{-4.3pt}\operatorname{I}}

\def\muc{\mu_{c}}

\def\nsize{{n}}

\def\rhor{\omega}

\def\La{\mathbb{L}}

\def\DP{D}
\def\DPc{\DP_{0}}
\def\DPb{\DP_{\rdb}}
\def\DPm{\DP_{\rdm}}

\def\gmi{\mathtt{b}}

\def\IF{\Bbb{F}}

\def\ex{\mathrm{e}}
\def\entrl{\mathbb{Q}}

\def\kullb{\cc{K}} 

\def\gm{\mathtt{g}}
\def\gmc{\gm_{c}}
\def\gmb{\gm}

\def\yy{\mathtt{y}}
\def\yyc{\yy_{c}}
\def\xx{\mathtt{x}}
\def\xxc{\xx_{c}}

\def\alp{\alpha}

\def\rups{\rr_{0}}

\def\Pdom{\mu_{0}}
\def\PDOM{\bb{\mu}_{0}}

\def\CS{\cc{E}}

\def\Ccb{m_{\rdb}}
\def\Ccm{m_{\rdm}}

\def\etas{\eta^{*}}

\def\nunu{\nu_{0}}

\def\dist{d}

\def\rdomega{\varrho}

\def\err{\diamondsuit}

\def\GP{G}

\def\LL{\cc{L}}

\def\dimp{p}

\def\BB{I\!\!B}
\def\vA{\mathtt{v}}

\def\spread{\Delta}

\def\thetav{\bb{\theta}}
\def\thetavs{\thetav^{*}}
\def\thetavc{\thetav'}
\def\thetavd{\thetav^{\circ}}

\def\thetas{\theta^{*}}

\def\thetavb{\thetav^{\dag}}

\def\vtheta{\vartheta}
\def\vthetav{\bb{\vtheta}}
\def\prior{\Pi}

\def\etavd{\etav^{\circ}}

\def\Sigmas{\Sigma_{0}}

\def\ups{\bb{\upsilon}}
\def\upss{\ups_{0}}

\def\Ups{\varUpsilon}
\def\Upsd{\Ups^{\circ}}
\def\Upss{\Ups_{\circ}}

\def\Thetas{\Theta_{0}}

\def\fs{f}

\def\eps{\epsilon}

\def\UU{\cc{Y}}

\def\UP{\cc{U}}

\def\VP{V}
\def\VPc{\VP_{0}}

\def\VV{H}

\def\VVc{\VV_{0}}

\def\vp{\mathbf{v}}

\def\lambdaB{{\lambda}^{*}}

\def\fis{\mathfrak{a}}

\def\B{\cc{B}}

\def\rr{\mathtt{r}}
\def\rrb{\rr^{*}}

\def\zz{\mathfrak{z}}

\def\score{\nabla}

\def\AssId{\mathcal{I}}

\def\dimn{\dimp_{\nsize}}

\def\CR{\mathcal{C}}

\def\vthetavb{\bar{\vthetav}}
\def\Covpost{\mathfrak{S}}

\def\Ec{\E^{\circ}}
\def\ff{f}

\def\entrlq{\entrl_{1}}

\def\kb{k^{*}}

\def\rderr{\chi}
\def\Excgr{\diamondsuit}

\def\xxn{\xx_{\nsize}}

\def\dimq{q}

\def\dimq{q}
\def\QQ{\mathbb{H}}
\def\QQg{\QQ_{2}}
\def\QQq{\QQ_{1}}

\def\lambdaB{\lambda_{\BB}}

\def\lambdav{\bb{\lambda}}
\def\etavd{\etav_{\circ}}
\def\thetavb{\breve{\thetav}}
\def\vthetavd{\Ec \vthetav}
\def\Covd{S_{\circ}}
\def\Covpostd{\Covpost_{\circ}}
\def\IS{\mathcal{I}}
\def\etas{\eta^{*}}
\def\Po{\operatorname{Po}}
\def\IF{\Bbb{F}}
\def\etavb{\bar{\etav}}
\def\etavd{\etav^{\circ}}
\def\Pc{\P^{\circ}}
\def\xxn{\xx_{\nsize}}
\def\CRd{\CR^{\circ}}

\def\CONST{\mathtt{C} \hspace{0.1em}}

\def\dimB{\mathtt{p}_{\BB}}

\def\nub{\nu}

\def\dimq{q}

\def\QQ{\mathbb{H}}
\def\QQg{\QQ_{2}}
\def\QQq{\QQ_{1}}

\def\qq{z}

\def\qqQ{\qq_{\QQ}}

\def\rderr{\chi}
\def\Excgr{\diamondsuit}

\def\Ccb{m}
\def\Ccm{m}
\def\BB{B}

\def\vthetavd{\vthetav^{\circ}}
\def\Indru{\Ind_{\rups}}

\def\BBh{U}
\def\betav{\bb{\beta}}
\def\DD{U}
\def\hsp{\tau}
\def\fiD{a}

\def\Prior{\Pi}
\def\prior{\pi}

\usepackage{amsmath,amssymb,amsthm}
\usepackage{natbib}
\usepackage{epsfig,graphicx}
\usepackage{comment}
\usepackage{color}
\usepackage{srcltx}
\usepackage[mathscr]{eucal}
\usepackage[math]{easyeqn}
\usepackage{etoolbox}

\ifims{
\textheight=23cm
\textwidth=14.8cm
\topmargin=0pt
\oddsidemargin=1.0cm
\evensidemargin=1.0cm
\linespread{1.3}
\renewenvironment{abstract}
    {\centerline{\textbf{Abstract}}\bigskip
      \begin{center}
       \begin{minipage}{11cm}
        \begin{small}
    }
    {   \end{small}
       \end{minipage}
      \end{center}
     \bigskip
    }

}{ 
}

\numberwithin{equation}{section}
\numberwithin{figure}{section}
\newcounter{example}[section]
\numberwithin{example}{section}
\newcounter{remark}[section]
\numberwithin{remark}{section}
\newtheorem{theorem}{Theorem}[section]
\newtheorem{proposition}[theorem]{Proposition}
\newtheorem{lemma}[theorem]{Lemma}
\newtheorem{corollary}[theorem]{Corollary}

\newtheorem{exmp}[example]{Example}
\newtheorem{rmrk}[remark]{Remark}
\newenvironment{example}{\begin{exmp}\rm}{\end{exmp}}
\newenvironment{remark}{\begin{rmrk}\rm}{\end{rmrk}}

\begin{document}
\thispagestyle{empty}
\ifims{
\title{\thetitle}
\ifau{ 
  \author{
    \authora
    \ifdef{\thanksa}{\thanks{\thanksa}}{}
    \\[5.pt]
    \addressa \\
    \texttt{ \emaila}
  }
}
{  
  \author{
    \authora
    \ifdef{\thanksa}{\thanks{\thanksa}}{}
    \\[5.pt]
    \addressa \\
    \texttt{ \emaila}
    \and
    \authorb
    \ifdef{\thanksb}{\thanks{\thanksb}}{}
    \\[5.pt]
    \addressb \\
    \texttt{ \emailb}
  }
}
{   
  \author{
    \authora
    \ifdef{\thanksa}{\thanks{\thanksa}}{}
    \\[5.pt]
    \addressa \\
    \texttt{ \emaila}
    \and
    \authorb
    \ifdef{\thanksb}{\thanks{\thanksb}}{}
    \\[5.pt]
    \addressb \\
    \texttt{ \emailb}
    \and
    \authorc
    \ifdef{\thanksc}{\thanks{\thanksc}}{}
    \\[5.pt]
    \addressc \\
    \texttt{ \emailc}
  }
}

\maketitle
\pagestyle{myheadings}
\markboth
 {\hfill \textsc{ \small \theruntitle} \hfill}
 {\hfill
 \textsc{ \small
 \ifau{\runauthora}
      {\runauthora and \runauthorb}
      {\runauthora, \runauthorb, and \runauthorc}
 }
 \hfill}
\begin{abstract}
\theabstract
\end{abstract}

\ifAMS
    {\par\noindent\emph{AMS 2000 Subject Classification:} Primary \kwdp. Secondary \kwds}
    {\par\noindent\emph{JEL codes}: \kwdp}

} 
{ 
\begin{frontmatter}
\title{\thetitle}


\runtitle{\theruntitle}

\ifau{ 
\begin{aug}
    \author{\authora\ead[label=e1]{\emaila}}
    \address{\addressa \\
     \printead{e1}}
\end{aug}

 \runauthor{\runauthora}
\affiliation{\affiliationa} }
{ 
\begin{aug}
    \author{\authora\ead[label=e1]{\emaila}\thanksref{t21}}
    \and
    \author{\authorb\ead[label=e2]{\emailb}\thanksref{t22}}
    
    \address{\addressa \\
     \printead{e1}}
    \address{\addressb \\
     \printead{e2}}
    \thankstext{t21}{\thanksa}
    \thankstext{t22}{\thanksb}
    \affiliation{\affiliationa, \affiliationb} 
    \runauthor{\runauthora and \runauthorb}
\end{aug}
} 
{ 
\begin{aug}
    \author{\authora\ead[label=e1]{\emaila}\thanksref{t21}}
    \and
    \author{\authorb\ead[label=e2]{\emailb}\thanksref{t22}}
    \and
    \author{\authorc\ead[label=e3]{\emailc}\thanksref{t23}}
    
    \address{\addressa \\
     \printead{e1}}
    \address{\addressb \\
     \printead{e2}}
    \address{\addressc \\
     \printead{e3}}
    \thankstext{t21}{\thanksa}
    \thankstext{t22}{\thanksb}
    \thankstext{t23}{\thanksc}
    \affiliation{\affiliationa, \affiliationb, \affiliationc} 
    \runauthor{\runauthora, \runauthorb, and \runauthorc}
\end{aug}}

\begin{abstract}
\theabstract
\end{abstract}



\end{frontmatter}
} 


\def\spreads{\spread_{\circ}}

\Chapter{Introduction}
\label{Swilksint}
The Fisher and Wilks Theorems are probably the most prominent results in parametric statistics.
The Fisher Theorem describes the asymptotic expansion of the maximum likelihood estimator (MLE) and 
in particular, claims its asymptotic normality.
The Wilks result describes the limiting behavior of the excess which is the maximum of the 
log-likelihood minus its value at the true point.
The limiting distribution is \( \chi^{2} \) with \( \dimp \) degrees of freedom,
where \( \dimp \) is the dimension of the parameter space.
These two results play a fundamental role in the whole theory of mathematical statistics and 
motivate or justify the most of statistical methods and approaches including
the estimation and testing theory, model selection, penalized procedures etc. 
%
%
Another prominent statistical result usually referred to as 
Bernstein -- von Mises (BvM) Theorem 
claims that the posterior measure is asymptotically normal
with the mean close to the maximum likelihood estimator (MLE) and the variance close to
the variance of the MLE.
This explains why this result is often considered as the Bayesian counterpart of the
frequentist Fisher Theorem about asymptotic normality of the MLE.
%
%
The BvM result provides a theoretical background for different Bayesian procedures.
In particularly, one can use Bayesian computations for evaluation of the MLE and its variance.
Also one can build elliptic credible sets 
using the first two moments of the posterior.

The classical versions of the Fisher, Wilks, and BvM Theorems are stated for the
standard parametric setup with a fixed parametric model and large samples.
We refer to \cite{VW1996}, \cite{vdV00} for a detailed historical overview.
The asymptotic results can be extended to a rather general statistical setup under the LAN condition 
by L. Le Cam; see, e.g., \cite{IH1981, LY2000}.
However, the theory is limited to the situation with a fixed parametric structure and large samples. 
Modern statistical applications force to consider statistical problems with a small or moderate 
sample size and large parameter dimension. 
In particular, a growing parameter dimension is naturally incorporated in the sieve approach, 
when the underlying nonparametric object is approximated by a sequence of parametric models with 
a large number of parameters;
see e.g. 
\cite{BiMa1993}, 
\cite{ShWo1994,shen1997}, 
\cite{vdG1993,vdG2002}.

An exact description of a complex problem is unrealistic even if (infinitely) many parameters 
are involved. 
This leads to another important issue -- model misspecification. 
It appears that the standard LAN arguments do not apply in this situation because 
all the small terms in the LAN expansion may depend on the dimension \( \dimp \)
and many steps require the true parametric structure.

The modern statistical theory faces the problem of making reliable statistical decisions 
in situations with small or moderate sample under possible model misspecification.
Recently \cite{SP2011} offered a version of the Le Cam LAN theory which applies to  
{finite samples}, {large parameter dimension}, and is robust against possible 
{model misspecification}. 
The proposed bracketing approach allows to describe explicitly all error terms in the 
expansion of the maximum likelihood estimator and of the corresponding excess.
The general results are illustrated by various examples including 
i.i.d. model, generalized linear and median regression.

This paper makes a further step in this direction by addressing two important questions:
the critical parameter dimension, and the use of the Fisher and Wilks expansion for
establishing the BvM Theorem.
%
The obtained bounds help to quantify and explicitly describe the error terms in all obtained results. 
In the situation when the parameter dimension grows with the sample size, one can use these 
results for describing the critical parameter dimension. 
In the important i.i.d. case, it appears that the MLE is consistent under 
``\( \dimp/n \) small'',  the Fisher expansion for the MLE valid under 
``\( \dimp^{2}/n \) small'', while the Wilks and BvM results apply under
``\( \dimp^{3}/n \) small''.
This founding differs from the existing results on dimensional asymptotics;
cf. \cite{Portnoy1984, Portnoy1985} for the M-estimation and \cite{Gh1999} for the BvM
result.
See Section~\ref{Scritdim} below for a more detailed comparison.
Further, we show that the BvM phenomenon only relies to the local quadratic approximation 
of the log-likelihood and some rough large deviation bounds. 
Any asymptotic arguments like the CLT or weak convergence are not required. 
This allows to include the case of a large parameter dimension and a possible model misspecification.
The main technical tool in the whole approach is a new bound on the maximum of 
a vector random field, which allows 
to improve and to simplify the bracketing device of \cite{SP2011}.

%

%

%
%

First we specify our set-up. 
Let \( \Yv \) denote the observed data and \( \P \) mean their distribution.
A general parametric assumption (PA) means that
\( \P \) belongs to \( \dimp \)-dimensional family
\( (\P_{\thetav}, \thetav \in \Theta \subseteq \R^{\dimp}) \) dominated by
a measure \( \PDOM \).
This family yields the log-likelihood function
\( L(\thetav) = L(\Yv,\thetav) \eqdef \log \frac{d\P_{\thetav}}{d\PDOM}(\Yv) \).
The PA can be misspecified, so, in general, \( L(\thetav) \) is a
\emph{quasi log-likelihood}.
The classical likelihood  principle suggests to estimate  \( \thetav \) by
maximizing the function \( L(\thetav) \):
\begin{EQA}[c]
    \tilde{\thetav}
    \eqdef
    \argmax_{\thetav \in \Theta} L(\thetav) .
\label{tthetamkGP}
\end{EQA}
If \( \P \not\in \bigl( \P_{\thetav} \bigr) \), then
the (quasi) MLE estimate \( \tilde{\thetav} \) from \eqref{tthetamkGP} is
still meaningful and it appears to be an estimate of the value \( \thetavs \)
defined by maximizing the expected value of \( L(\thetav) \):
\begin{EQA}[c]
    \thetavs
    \eqdef
    \argmax_{\thetav \in \Theta} \E L(\thetav).
\label{thetavsd}
\end{EQA}
Here
\( \thetavs \) is the true value in the parametric situation and can be viewed as the
parameter of the best parametric fit in the general case.
The study is non-asymptotic, that is, we proceed with only one sample \( \Yv \).
One can easily extend it to an asymptotic setup in which the data, its distribution,
the parameter space and the parametric family depend on the asymptotic parameter like 
the sample size. 
One example is given below in Section~\ref{Scritdim} for the case of an i.i.d. sample.

The Fisher expansion of the qMLE \( \tilde{\thetav} \) is given as follows:
\begin{EQA}[c]
    \DPc (\tilde{\thetav} - \thetavs)
    \approx
    \xiv
    \eqdef
    \DPc^{-1} \nabla L(\thetavs),
\label{DPcttxiv}
\end{EQA}
where \( \nabla L(\thetav) = \frac{dL}{d\thetav}(\thetav) \) and
\( \DPc^{2} \eqdef - \nabla^{2} \E L(\thetavs) \) is the analog of the total
Fisher information matrix.
In classical situations, the standardized score \( \xiv \) is asymptotically
standard normal yielding asymptotic root-n normality and efficiency of the MLE
\( \tilde{\thetav} \).
Theorem~\ref{Tconflocro} carefully describes how the error of this expansion depends 
on the parameter dimension \( \dimp \) and the regularity of the model.
%
The Wilks expansion means 
\begin{EQA}
	L(\tilde{\thetav}) - L(\thetavs)
	& \approx &
	\| \xiv \|^{2}/2 .
\label{LttLts}
\end{EQA}
Again, if the vector \( \xiv \) is asymptotically standard normal, 
the expansion yields the classical \( \chi^{2}_{\dimp} \) asymptotic distribution 
for the excess \( L(\tilde{\thetav}) - L(\thetavs) \).


In the Bayes setup \( \vthetav \) is a random element following a prior
measure \( \Prior \) on the parameter set \( \Theta \).
The posterior describes the conditional distribution of \( \vthetav \) given \( \Yv \)
obtained by normalization of the product \( \exp \bigl\{ L(\thetav) \bigr\} \Prior(d\thetav) \).
This relation is usually written as
\begin{EQA}
    \vthetav \cond \Yv
    & \propto &
    \exp \bigl\{ L(\thetav) \bigr\} \, \Prior(d\thetav) .
\label{Bayesform}
\end{EQA}
A model misspecification means that \( L(\thetav) \) is not necessarily a true log-likelihood.
In this case, the Bayes formula \eqref{Bayesform} describes a
\emph{quasi posterior}; cf. \cite{Ch2003}.
Section~\ref{SSnoninf} studies some general properties of the posterior measure
focusing on the case of a \emph{non-informative} prior \( \Prior \).
The main result claims that the distribution of
\( \DPc (\vthetav - \thetavs) - \xiv \) given \( \Yv \) is nearly standard normal.
Comparing with \eqref{DPcttxiv} indicates that the posterior is nearly centered at
the qMLE \( \tilde{\thetav} \) and its variance is close to \( \DPc^{-2} \).
So, our result extends the BvM theorem to the considered general setup.
Section~\ref{Sregularprior} comments how the result can be extended to the case of any 
\emph{regular} prior.

The whole study is nonasymptotic and all ``small'' terms are carefully described.
This helps to understand how the parameter dimension is involved and particularly
to address the question of a \emph{critical dimension}; see
Section~\ref{SBvMiid} which specifies the result to the i.i.d. case with \( \nsize \) observations
and links the obtained results to the classical literature. 

The paper is structured as follows:
Section~\ref{Scondgllo} presents 
the imposed condition, the frequentist results are given in Section~\ref{Sfreqres}, 
the BvM results and its extensions are described in Section~\ref{SGaussapprox}.
Section~\ref{SauxBvM} presents some auxiliary results 
which can be of independent interest 
and provides the proofs.
Some useful auxiliary facts are collected in the Appendix.

\Chapter{Main results}
\label{Smainresults}
This section presents our main results which include the Fisher and Wilks expansions and 
the Bernstein--von Mises Theorem for 
a non-classical and non-asymptotic framework.
First we present the frequentist results: concentration and large deviation 
properties of the maximum likelihood estimator \( \tilde{\thetav} \),
the Fisher expansion for the difference \( \tilde{\thetav} - \thetav \) and the Wilks expansion 
for the excess \( L(\tilde{\thetav}) - L(\thetav) \).
Then we switch to the Bayesian framework and establish the BvM result for the special case 
of a non-informative prior
given by the Lebesgue measure \( \prior(\thetav) \equiv 1 \) on \( \R^{\dimp} \)
and extend the results to any prior with a continuous density on a vicinity
of the central point \( \thetavs \).
The results are stated in a concise way, all the terms are given explicitly. 
Surprisingly, the leading terms in all bounds are sharp, 
in particular, the classical results on asymptotic efficiency can be easily derived from the obtained 
expansions. 

Introduce the notation \( L(\thetav,\thetavs) = L(\thetav) - L(\thetavs) \) for
the (quasi) log-likelihood ratio.
The main step in the approach 
is the following \emph{uniform local bracketing result}:
\begin{EQA}[c]
    \La(\thetav,\thetavs) - \spread
    \le
    L(\thetav,\thetavs)
    \le
    \La(\thetav,\thetavs) + \spread,
    \qquad
    \thetav \in \Thetas .
\label{LLLmbin}
\end{EQA}
Here \( \La(\thetav,\thetavs) \) is a quadratic in
\( \thetav - \thetavs \) expression,
\( \spread \) is a small error only depending on \( \Thetas \) which is a local vicinity
of the central point \( \thetavs \).
This result can be viewed as an extension of the famous Le Cam
\emph{local asymptotic normality} (LAN) condition.
The LAN condition postulates an approximation of the
log-likelihood \( L(\thetav) \) by a nearly Gaussian process; see e.g. \cite{IH1981}
or \cite{KlVa2012} for an extension of this condition (stochastic LAN).
The bracketing bound \eqref{LLLmbin} requires only   
some general conditions listed in Section~\ref{Scondgllo}.
A model misspecification case is included.
Similarly to the LAN theory, the bracketing result has a number of remarkable
corollaries like the Wilks and Fisher Theorems; see Theorems~\ref{Tconflocro} and 
\ref{TWilks2r}.
Below we show that the BvM result is in some sense also a corollary of \eqref{LLLmbin}.

For making a precise statement, we have to specify the ingredients of the bracketing
device.
The most important one is a symmetric positive \( \dimp \times \dimp \)-matrix
\( \DPc^{2} \). 
In typical situations, it can be defined as the negative Hessian of the expected 
log-likelihood: \( \DPc^{2} = - \nabla^{2} \E L(\thetavs) \).
Also one has to specify a radius \( \rups \) entering
in the definition of the local vicinity
\( \Thetas(\rups) \) of the central point \( \thetavs \):
\( \Thetas(\rups) = \{ \thetav \colon \| \DPc (\thetav - \thetavs) \| \le \rups \} \).
The bracketing result \eqref{LLLmbin} can be stated for \( \Thetas = \Thetas(\rups) \) 
with
\begin{EQA}
    \La(\thetav,\thetavs)
    & \eqdef &
    (\thetav - \thetavs)^{\T} \nabla L(\thetavs)
    - \| \DPc (\thetav - \thetavs) \|^{2}/2
    \\
    &=&
    \xiv^{\T} \DPc (\thetav - \thetavs)
    - \| \DPc (\thetav - \thetavs) \|^{2}/2
\label{bLquadlocpr}
\end{EQA}
and
\begin{EQA}
    \xiv
    & \eqdef &
    \DPc^{-1} \nabla L(\thetavs) .
\label{xivalpspr}
\end{EQA}
The construction is essentially changed relative to \cite{SP2011} 
(in fact, it is simplified) by using only one matrix 
\( \DPc^{2} \) while \cite{SP2011} used three matrices \( \DPb^{2} \), \( \DPm^{2} \), and 
\( \VPc^{2} \).
%
The bracketing bound \eqref{LLLmbin} becomes useful if the error \( \spread \) is 
relatively small and can be neglected.

\Section{Conditions} 
\label{Scondgllo} 
This section collects the conditions which
are systematically used in the text. 
The conditions are quite general and seem to be non-restrictive; see the discussion 
at the end of the section.
We mainly require some regularity and smoothness of the 
log-likelihood process \( L(\thetav) \).
With \( \DPc^{2} = - \nabla^{2} \E L(\thetavs) \), define the local elliptic sets 
\( \Thetas(\rr) \) as
\begin{EQA}[c]
    \Thetas(\rr)
    \eqdef
    \{ \thetav \in \Theta \colon \| \DPc (\thetav - \thetavs) \| \le \rr \} .
\label{Theta0R}
\end{EQA}
We distinguish between local and global conditions.
Local ones are stated on \( \Thetas(\rups) \), while the global one corresponds to 
\( \rr \geq \rups \), where the value \( \rups \)
will be specified later.

The first condition requires that
the expected log-likelihood \( \E L(\thetav) \) is twice continuously differentiable.

\begin{description}
    \item[\( \bb{(\LL_{0})} \)]
    \textit{
    For each \( \rr \leq \rups \), 
    there is a constant \( \rddelta(\rr) \leq 1/2 \) such that
    it holds for any \( \thetav \in \Thetas(\rr) \) and 
    \( \DP^{2}(\thetav) \eqdef - \nabla^{2} \E L(\thetav) \):
    }
\begin{EQA}[c]
\label{LmgfquadELGP}
    \bigl\|
		\DPc^{-1} \DP^{2}(\thetav) \DPc^{-1} - \Id_{\dimp} 
    \bigr\|_{\infty}
    \le
    \rddelta(\rr)  .
\end{EQA}
\end{description}
Under \( (\LL_{0}) \), it follows from the second order Taylor
expansion at \( \thetavs \):
\begin{EQA}
    \bigl|
        - 2 \E L(\thetav,\thetavs)
        - \| \DPc (\thetav - \thetavs) \|^{2}
    \bigr|
    & \le &
    \rddelta(\rr) \| \DPc (\thetav - \thetavs) \|^{2} ,
    \quad
    \thetav \in \Thetas(\rr).
\label{EdeltGP}
\end{EQA}

For \( \rr > \rups \),  we need a global identification property which ensures that the deterministic component 
\( \E L(\thetav,\thetavs) \) of the log-likelihood is competitive with the
variation of the stochastic component.

\begin{description}
    \item[\( \bb{(\LL)} \)]
    \textit{
    There exists \( \gmi(\rr) > 0 \) 
    such that
    \( \rr \gmi(\rr) \) is non-decreasing and
	}
\begin{EQA}
    \frac{- 2 \E L(\thetav,\thetavs)}{\| \DPc (\thetav - \thetavs) \|^{2}}
    & \ge &
    \gmi(\rr) ,
    \qquad
    \forall \rr \geq \rups, \, \, \thetav \in \Thetas(\rr) .
\label{xxentrttGP}
\end{EQA}    
\end{description}

Now we consider the stochastic component of the process \( L(\thetav) \):
\begin{EQA}[c]
	\zeta(\thetav)
	\eqdef
	L(\thetav) - \E L(\thetav) .
\end{EQA}
We assume that it is twice differentiable and denote by \( \nabla \zeta(\thetav) \) its gradient and by 
\( \nabla^{2} \zeta(\theta) \) its Hessian matrix. 

\begin{description}
\item[\( \bb{(E\!D_{0})} \)]
    \emph{ There exist a positive symmetric matrix \( \VPc^{2} \),
    and constants \( \gmb > 0 \), \( \nunu \ge 1 \) such that
    \( \Var\bigl\{ \nabla \zeta(\thetavs) \bigr\} \le \VPc^{2} \) and 
    }
\begin{EQA}[c]
\label{expzetaclocGP}
	\sup_{\gammav \in \R^{\dimp}} 
    \log \E \exp\biggl\{
        \lambda \frac{\gammav^{\T} \nabla \zeta(\thetavs)}
                     {\| \VPc \gammav \|}
    \biggr\} \le
    \nunu^{2} \lambda^{2} / 2, \qquad 
    |\lambda| \le \gmb .
\end{EQA}

\item[\( \bb{(E\!D_{2})} \)]
    \emph{There exist a value \( \rhor > 0 \) 
    and for each \( \rr > 0 \), a constant \( \gm(\rr) > 0 \) such that
    it holds for any \( \thetav \in \Thetas(\rr) \):    }
\begin{EQA}[c]
    \sup_{\gammav_{1},\gammav_{2} \in \R^{\dimp} }
    \log \E \exp\biggl\{ 
    	\frac{\lambda}{\rhor} \,\,
        \frac{\gammav_{1}^{\T} \nabla^{2} \zeta(\thetav) \gammav_{2}} 
        	 {\| \DPc \gammav_{1} \| \cdot \| \DPc \gammav_{2} \|}
	\biggr\} 
    \le 
    \frac{\nunu^{2} \lambda^{2}}{2} \, ,
    \qquad 
    |\lambda| \leq \gm(\rr).
\label{gUUgem2}
\end{EQA}
\end{description}
Below we only need that  the constant \( \gm(\rr) \) 
is larger than \( \CONST \dimp \) for a fixed constant \( \CONST \) and all \( \rr \).
%

%

The \emph{identifiability condition} relates the matrices \( \VPc^{2} \) and \( \DPc^{2} \).
\begin{description}
  \item[\( (\bb{\AssId}) \)] 
      There is a constant 
      \( \fis > 0 \) such that 
\begin{EQA}
	\fis^{2} \DPc^{2} 
	& \ge & 
	\VPc^{2} . 
\label{lamGPDPVPfis}
\end{EQA}
  
\end{description}

\begin{remark}
The conditions involve some constants.
We distinguish between important constants and technical ones.
The impact of the important constants is shown in our results, the list includes 
\( \rddelta(\rr) \), \( \rhor \), and \( \fis \).
The constant \( \fis \) can be viewed as the largest eigenvalue of 
\( \BB = \DPc^{-1} \VPc^{2} \DPc^{-1} \) and it enters in the definition of 
the upper quantile function \( \qq(\BB,\xx) \) for \( \| \xiv \| \); 
see Proposition~\ref{LLbrevelocro} below.
The other constants like \( \nunu \) or \( \gm(\rr) \) are technical.
The constant \( \nunu \) is introduced for convenience only, it can be omitted by 
rescaling the matrix \( \VPc \). 
In the asymptotic setup it can usually be selected very close to one. 
\end{remark}

\begin{remark}
We  briefly comment how restrictive the imposed conditions are.
\cite{SP2011}, Section 5.1, considered in details the i.i.d. case and presented some mild 
sufficient conditions on the parametric family which imply the above general conditions.
Condition \( (E\!D_{0}) \) requires some exponential moments of the observations 
(errors). 
Usually one only assumes some finite moments of the errors; cf. \cite{IH1981}, Chapter~2.
Our condition is a bit more restrictive but it allows to obtain some finite sample 
bounds. 
Condition \( (\LL_{0}) \) only requires some regularity of the considered parametric 
family and is not restrictive.
Conditions \( (E\!D_{2}) \) with \( \gm(\rr) \equiv \gm > 0 \) and \( (\cc{L}) \) 
with \( \gmi(\rr) \equiv \gmi > 0 \)
are easy to verify if the parameter set 
\( \Theta \) is compact and the sample size \( \nsize \) exceeds \( \CONST \dimp \) 
for a fixed constant \( \CONST \).
It suffices to check a usual identifiability condition that the value 
\( \E L(\thetav,\thetavs) \) does not vanish for \( \thetav \neq \thetavs \).

The regression and generalized regression models are included as well; cf. 
\cite{Gh1999,Gh2000} or \cite{Kim2006}.
\cite{SP2011}, Section 5.2, argued that \( (E\!D_{2}) \) is automatically fulfilled for 
a generalized linear model, while \( (E\!D_{0}) \) requires that 
regression errors have to fulfill some exponential moments conditions. 
If this condition is too restrictive and a more stable (robust) estimation procedure is 
desirable, one can apply the LAD-type contrast leading to median regression. 
\cite{SP2011}, Section 5.3, showed for the case of linear median regression that 
all the required conditions are fulfilled automatically if the sample size \( \nsize \) 
exceeds \( \CONST \dimp \) for a fixed constant \( \CONST \).
\cite{SpWe2012} applied this approach for local polynomial quantile regression.
\cite{zaitsev2013properties} applied the approach to the problem of regression with Gaussian 
process where the unknown parameters enter in the likelihood in a rather complicated 
way. 
%
\end{remark}

\Section{Properties of the MLE \( \tilde{\thetav} \)}
\label{Sfreqres}
This section collects the main results about the MLE \( \tilde{\thetav} \).
We begin by a large deviation bound which ensures a small probability of the event
\( \tilde{\thetav} \not\in \Thetas(\rups) \).
Then we present the Fisher and Wilks expansions.
The formulation involves two growing functions of the argument \( \xx \):
\( \qq(\BB,\xx) \) and 
\( \qqQ(\xx) \).
The functions are given analytically and only depend on the parameters of the model.
The function \( \qq(\BB,\xx) \) describes the quantiles of the norm of the vector \( \| \xiv \| \) from 
\eqref{xivalpspr}.
The definition is given in \eqref{zzxxppdBlro}.
Further, the function \( \qqQ(\xx) \) is related to the entropy of the parameter space and it is
given by \eqref{zzxxgfindef}.
In typical situations one can use the upper bound \( \qq^{2}(\xx) \leq \CONST (\dimp + \xx) \) 
for both functions.
The first result explains the choice of \( \rups \) ensuring with a high probability that  
\( \tilde{\thetav} \in \Thetas(\rups) \).

\begin{theorem}
\label{TMLE}
Suppose \( (E\!D_{0}) \) and \( (E\!D_{2}) \), 
\( (\LL_{0}) \), \( (\LL) \), and \( (\AssId) \). 
Let also 
the function  \( \gmi(\rr) \) in \( (\LL) \) satisfy 
\begin{EQA}
    \gmi(\rr) \, \rr
    & \ge &
    2 \qq(\BB,\xx) + 2 \rdomega(\rr,\xx), 
    \quad
    \rr > \rups,
\label{cgmibrr}
\end{EQA}
where 
\begin{EQA}[c]
    \rdomega(\rr,\xx)
    \eqdef
    6 \nunu \, \qqQ\bigl(\xx + \log(2\rr/\rups) \bigr) \, \rhor  .
\label{Exceqrrrhdef}
\end{EQA}  
Then 
\begin{EQA}[c]
    \P\bigl( \tilde{\thetav} \not\in \Thetas(\rups) \bigr)
    \le 
    3 \ex^{-\xx} .
\label{PCAxxglrrsGP}
\end{EQA}    
\end{theorem}

\begin{remark}
The radius \( \rups \) has to fulfill 
\eqref{cgmibrr}.
This condition is easy to check in typical situations.
One can use that \( \gmi(\rups) \geq 1 - \rddelta(\rups) \approx 1 \),
that the constant \( \rhor \) is small, and
\( \rr \gmi(\rr) \) grows with \( \rr \).
A simple rule \( \rups \geq (2 + \delta) \qq(\BB,\xx) \) for some \( \delta > 0 \) 
works in most of cases. 
\end{remark}

%

Now we state the result about the Fisher expansion for the qMLE 
\( \tilde{\thetav} \). 

\begin{theorem}
\label{Tconflocro}
Suppose the conditions of Theorem~\ref{TMLE}.
On a random set 
\( \Omega(\xx) \) of a dominating probability at least \( 1 - 4 \ex^{-\xx} \), it holds 
\begin{EQA}[c]
    \| \DPc (\tilde{\thetav} - \thetavs) - \xiv \|
	\leq 
    \Excgr(\rups,\xx) ,
\label{ttusmxivF}
\end{EQA}
where for the function \( \qqQ(\xx) \) given by \eqref{zzxxgfindef}, 
the value \( \Excgr(\rr,\xx) \) is given by 
\begin{EQA}[c]
    \Excgr(\rr,\xx)
    \eqdef
    \bigl\{ \rddelta(\rr) 
    + 6 \nunu \, \qqQ(\xx) \, \rhor \bigr\}\, \rr .
\label{Exceqrrrxxhdef}
\end{EQA}  
%
\end{theorem}

Our version of the Wilks result can be stated in the following form.

\begin{theorem}
\label{TWilks2r}
Suppose the conditions of Theorem~\ref{TMLE}.
On a random set 
\( \Omega(\xx) \) of a dominating probability at least \( 1 - 5 \ex^{-\xx} \), it holds with 
\( \Excgr(\rups,\xx) \) from \eqref{Exceqrrrxxhdef}
\begin{EQA}
\label{WilLLe}
    \bigl| 2 L(\tilde{\thetav},\thetavs) - \| \xiv \|^{2} \bigr|
    & \le &
    2 \spread(\rups,\xx) ,
    \\
    \Bigl| 
    	\sqrt{ 2L(\tilde{\thetav},\thetavs) } 
		- \| \xiv \| 
	\Bigr|
    & \le &
    3 \Excgr(\rups,\xx) ,
\label{sqWilLL}
\end{EQA}    
where
\begin{EQA}
	\spread(\rups,\xx)
	& \eqdef &
    \bigl\{ 2 \rups + \qq(\BB,\xx) \bigr\} \, \Excgr(\rups,\xx) .
\label{spreadrxxdef}
\end{EQA}
\end{theorem}

\begin{remark}
The classical Fisher and Wilks results describe asymptotic behavior of the MLE 
\( \tilde{\thetav} \) and of the excess \( L(\tilde{\thetav},\thetavs) \).
The whole derivations are based on expansions similar to \eqref{ttusmxivF} and \eqref{WilLLe}
and on the limiting behavior of the vector \( \xiv \) and its squared norm.
Under standard assumptions in the regression or i.i.d. setup the vector \( \xiv \) is 
standard normal and \( \| \xiv \|^{2} \) is asymptotically \( \chi^{2} \) with \( \dimp \) degrees 
of freedom. 
The asymptotic distribution of the MLE \( \tilde{\thetav} \) or of the excess 
\( L(\tilde{\thetav},\thetavs) \) can be used for building the confidence sets or for test critical values. 
However, the use of asymptotic arguments is limited and faces serious problems in practical applications.

This especially concerns the likelihood ratio statistic \( L(\tilde{\thetav},\thetavs) \).
It is well recognized that the accuracy of \( \chi^{2} \)-approximation of the tails 
of \( L(\tilde{\thetav},\thetavs) \) is very poor and a reasonable quality requires a huge sample size. 
If the parameter dimension grows with \( \nsize \) this problem becomes even more crucial. 
The qualitative tail behavior of \( \| \xiv \|^{2} \) is described in Proposition~\ref{LLbrevelocro}
but the upper bounds given there appear to be too conservative for practical use. 
\end{remark}

\begin{remark}
Another issue is a possible model misspecification. 
The expansion \eqref{WilLLe} applies, 
even if \( L(\thetav) \) is a quasi-log-likelihood function.
However, the covariance matrix 
\( \VPc^{2} = \Var\bigl\{ \nabla L(\thetavs) \bigr\} \) of the score does not necessarily 
coincide with the information matrix \( \DPc^{2} \).
Then the covariance matrix of the vector \( \xiv \) follows the famous ``sandwich'' formula 
\( \Var(\xiv) = \DPc^{-1} \VPc^{2} \DPc^{-1} \), 
and the distribution of the squared norm \( \| \xiv \|^{2} \)
depends on the unknown covariance matrix \( \VPc^{2} \). 
\end{remark}

The results presented above focus on the expansions of the MLE \( \tilde{\thetav} \) and 
on the excess \( L(\tilde{\thetav},\thetavs) \).
Numerical results (not presented here) indicate that the accuracy of the expansions \eqref{ttusmxivF} and \eqref{WilLLe}
is very reasonable even for moderate and small samples and it is stable w.r.t. possible model misspecifications. 
It seems that such expansions are of independent interest and can be used for many 
further statistical tasks. 
The next section explains how the obtained results can be used for proving the prominent 
Bernstein -- von Mises Theorem. 
Similarly, one can use this technique 
for showing the consistency of a bootstrap resampling procedure.


\Section{A Bernstein - von Mises Theorem}
\label{SGaussapprox}
\label{SSnoninf}
This section discusses the properties of the posterior distribution for a non-informative 
and a regular prior and the log-likelihood function \( L(\thetav) \). 
The Bernstein -- von Mises result claims that the posterior is nearly normal
with the mean \( \tilde{\thetav} \) and the variance \( \DPc^{-2} \),
where \( \DPc^{2} \) is the information matrix from condition \( (\LL_{0}) \).
We refer to \cite{KlVa2012} for a detailed historical overview around the BvM result.
%
There is a number of papers in this direction recently appeared.
We mention
\cite{GhGhVa2000,GhVa2007} for a general theory in the i.i.d. case;
\cite{Gh1999},
\cite{Gh2000} for high dimensional linear models;
\cite{BoGa2009},
\cite{Kim2006} for some special non-Gaussian models;
\cite{Sh2002},
\cite{BiKl2012}, \cite{RiRo2012}, \cite{Ca2012}, \cite{CaRo2013} for 
a semiparametric version of the BvM result for different models;
\cite{KlVa2006},
\cite{BuMi1998}, for the misspecified parametric case,
\cite{CaRo2013},
among many others.
A general framework for the BvM result is given in \cite{KlVa2012} in terms of the so called stochastic LAN conditions.
This condition extends the classical LAN condition and basically means a kind of quadratic expansion of the log-likelihood in a root-n vicinity of the central point. 
The approach applies even if the parametric assumption is misspecified, however, it requires a fixed parametric model and large samples. 
Extensions to nonparametric models with infinite or growing 
parameter dimension \( \dimp \) exist for some special situations, see e.g. \cite{Fr1999} and 
\cite{Gh1999,Gh2000} for linear models or \cite{Bo2011} for Gaussian regression,
or \cite{CaNi2013} for the functional Gaussian case.
Though the main arguments behind BvM results are similar in all studies, the way of bounding 
the error terms in the BvM results are essentially different.
The approach of this paper allows to get explicit upper bounds on the error of Gaussian 
approximation for the posterior law which apply for finite samples and 
admit model misspecification.
In a special case of an i.i.d. sample, one can precisely control how the error terms depend on the sample size \( \nsize \) and the dimension \( \dimp \) and judge about the applicability 
of the approach when \( \dimp \) grows with \( \nsize \).
We also show that the posterior mean is a very good approximation of the MLE,
while the posterior variance estimates the inverse of the Fisher information matrix.
Section~\ref{SellCS} discusses a possible construction of credible sets based 
on the posterior mean and variance.
Our results are stated for the non-informative prior.
In the Bayesian nonparametric literature the contraction rate is heavily influenced by the prior.
However, in the considered setup, the prior structure does not significantly affect the results 
and the main statements continue to hold for any a regular prior; see Section~\ref{Sregularprior}.

Let \( \vthetav \) mean a random element on the parameter set \( \Theta \),
by \( \prior(\thetav) \) we denote a prior density.
In this section we assume that \( \vthetav \) is uniformly distributed on \( \Theta \) 
with \( \prior(\thetav) \equiv 1 \).
The posterior distribution of \( \vthetav \) given \( \Yv \) is described by the product 
density \( \exp\bigl\{ L(\thetav) \bigr\} \) normalized by the marginal density 
\( p(\Yv) = \int_{\Theta} \exp\bigl\{ L(\thetav) \bigr\} \prior(\thetav) d\thetav \).
Introduce the posterior moments
\begin{EQA}[c]
    \vthetavb
    \eqdef
    \E \bigl( \vthetav \cond \Yv \bigr) ,
    \quad
    \Covpost^{2}
    \eqdef
    \Cov (\vthetav \cond \Yv)
    \eqdef
    \E \bigl\{ (\vthetav - \vthetavb) (\vthetav - \vthetavb)^{\T} \cond \Yv \bigr\}.
    \qquad
\label{vtheCovvth}
\end{EQA}
An important feature of the posterior distribution is that it 
can be numerically assessed.
In particular, one can evaluate its moments from \eqref{vtheCovvth}.
If we know in addition that the posterior is nearly normal, 
then the posterior is completely specified.
This information can be effectively used for building Bayesian credible sets 
with an elliptic shape; see the next section.

Before stating the results, introduce some more notations.
Define 
\begin{EQA}
	\thetavb 
	&=&
	\thetavs + \DPc^{-2} \score L(\thetavs) 
	=
	\thetavs + \DPc^{-1} \xiv.
\label{thetavdBvM}
\end{EQA}
The result of Theorem~\ref{Tconflocro} implies the expansion of the MLE
\( \tilde{\thetav} \) in the form \( \| \DPc (\tilde{\thetav} - \thetavb) \| \leq \err(\rups,\xx) \).
%
This section presents a version of the BvM result in the considered nonasymptotic setup
which claims that \( \vthetavb \) is close to \( \thetavb \) and thus to \( \tilde{\thetav} \), 
\( \Covpost^{2} \) is nearly equal to \( \DPc^{-2} \),
and \( \DPc \bigl( \vthetav - \thetavb \bigr) \) is nearly standard
normal conditionally on \( \Yv \).

\begin{theorem}
\label{TBvM}
Suppose the conditions of Theorem~\ref{TMLE}.
Let also \( \gmi(\rr) \) from \( (\LL) \) satisfies
\begin{EQA}
	\rr^{2} \gmi^{2}(\rr)
	& \geq &
	\xx
	+
	2 \dimp 
	+ 4 \qq^{2}(\BB,\xx) 
	+ 8 \rr \, \gmi(\rr) \rdomega(\rr,\xx),
	\qquad
	\rr \geq \rups,
\label{rrrupsboundco}
\end{EQA}
with \( \rdomega(\rr,\xx) \) from \eqref{Exceqrrrhdef}.
Then it holds on a random set \( \Omega(\xx) \) of probability at least \( 1 - 5 \ex^{-\xx} \)
\begin{EQA}
	\P\bigl( \vthetav \not\in \Thetas(\rups) \cond \Yv \bigr)
	& \leq &
	\ex^{-\xx}.
\label{contracrate}
\end{EQA}
Also on \( \Omega(\xx) \) with \( \spreads(\xx) = \rups \Excgr(\rups,\xx) \),
cf. \eqref{spreadrxxdef};
\begin{EQA}
\label{vthetabtherdb}
    \| \DPc (\vthetavb - \thetavb) \|^{2}
    & \le &
    4 \spreads(\xx) + 4 \ex^{-\xx},
    \\
    \bigl\| \Id_{\dimp} - \DPc \Covpost^{2} \DPc \bigr\|_{\infty}
    & \le &
    4 \spreads(\xx) + 4 \ex^{-\xx} .
\label{CovpostBvM}
\end{EQA}
Moreover, for any \( \lambdav \in \R^{\dimp} \) with \( \| \lambdav \|^{2} \le \dimp \)
\begin{EQA}[c]
    \Bigl|
        \log \E \Bigl[
            \exp\bigl\{ \lambdav^{\T} \DPc (\vthetav - \thetavb) \bigr\}
            \cond \Yv
        \Bigr]
        - \| \lambdav \|^{2}/2
    \Bigr|
    \le
    2 \spreads(\xx) + \ex^{-\xx} ,
\label{BvMmainr}
\end{EQA}
and for any measurable set \( A \subset \R^{\dimp} \)
\begin{EQ}[rcr]
	\P\bigl( \DPc (\vthetav - \thetavb) \in A \cond \Yv \bigr)
	& \geq &
	\exp \bigl\{ - 2 \spreads(\xx) - 3 \ex^{-\xx} \bigr\}
	\P\bigl( \gammav \in A \bigr) - \ex^{- \xx},
	\\
	\P\bigl( \DPc (\vthetav - \thetavb) \in A \cond \Yv \bigr)
	& \leq &
	\exp \bigl\{ 2 \spreads(\xx) + 2 \ex^{-\xx} \bigr\}
	\P\bigl( \gammav \in A \bigr) + \ex^{- \xx}.
\label{PcPgamGAm}
\end{EQ}
\end{theorem}

One can see that all statements of Theorem~\ref{TBvM} require 
``\( \spreads(\xx) = \rups \Excgr(\rups,\xx) \) 
is small''. 
Later we show that the results continue to hold if \( \thetavb \) 
is replaced by any efficient estimate \( \hat{\thetav} \), 
e.g. by the MLE \( \tilde{\thetav} \), satisfying 
\( \| \DPc (\hat{\thetav} - \thetavb) \| \leq \rups \) 
with a dominating probability.

\begin{remark}
The BvM result is stated under essentially the same list of conditions 
as the frequentist results of Theorem~\ref{TMLE} through \ref{TWilks2r}.
Similarly to the previous results, 
the normal approximation of the posterior is entirely based on the smoothness properties 
of the likelihood function and does not involve any asymptotic arguments like 
weak convergence or convergence in probability, or the Central Limit Theorem.
\end{remark}

\begin{remark}
The bound \eqref{rrrupsboundco} is very similar to the bound \eqref{cgmibrr} of Theorem~\ref{TMLE}.
It can be spelled out as the conditions that \( \rups^{2} \geq {2 \dimp + \xx} + 4 \zz^{2}(\BB,\xx) \),
\( \gmi(\rups) \approx 1 \), and \( \rr \gmi(\rr) \) grows with \( \rr \).
\end{remark}

\Subsection{The use of posterior mean and variance for credible sets}
\label{SellCS}

This section discusses a possibility of building some Bayesian credible sets in the 
elliptic form motivated by the Gaussian approximation of the posterior. 
The BvM result ensures that the posterior can be well approximated by the normal law with the mean 
\( \thetavb \) and the covariance \( \DPc^{-2} \).
%
This means that the posterior probability of the set 
\begin{EQA}[c]
    \CRd(A)
    =
    \bigl\{ \thetav \colon \DPc (\thetav - \thetavb) \in A \bigr\} 
\label{credsetBa}
\end{EQA}
is close to \( \P(\gammav \in A) \) up to an error term of order \( \spreads(\xx) \).
Unfortunately, the quantities \( \thetavb \) and \( \DPc^{2} \) are unknown and cannot be used for building the elliptic credible sets. 
A natural question is whether one can replace these values by some empirical counterparts without any substantial change of the posterior mass. 
An answer is given by the following result.

\begin{theorem}
\label{TBvMhat}
Let a vector \( \hat{\thetav} \) and a symmetric matrix \( \hat{\DP}^{2} \) fulfill
\begin{EQA}
	\| \DPc (\hat{\thetav} - \thetavb) \|	
	& \leq &
	\beta,
	\qquad
	\hat{\DP}^{2}
	\leq
	\fiD^{2} \DPc^{2}, 
	\qquad
	\tr \bigl( \DPc^{-1} \hat{\DP}^{2} \DPc^{-1} - \Id_{\dimp} \bigr)^{2}
	\leq 
	\eps^{2} .
\label{DPchatDPbeta}
\end{EQA}
Then with \( \hsp = \frac{1}{2} \sqrt{\fiD^{2} \beta^{2} + \eps^{2} } \), it holds on a set
\( \Omega(\xx) \) of probability \( 1 - 5 \ex^{-\xx} \)
\begin{EQ}[rcr]
	\P\bigl( \hat{\DP} (\vthetav - \hat{\thetav}) \in A \cond \Yv \bigr)
	& \geq &
	\exp \bigl( - 2 \spreads(\xx) - 3 \ex^{-\xx} \bigr)
	\bigl\{ \P\bigl( \gammav \in A \bigr) - \hsp \bigr\} - \ex^{- \xx},
	\\
	\P\bigl( \hat{\DP} (\vthetav - \hat{\thetav}) \in A \cond \Yv \bigr)
	& \leq &
	\exp \bigl( 2 \spreads(\xx) + 2 \ex^{-\xx} \bigr)
	\bigl\{ \P\bigl( \gammav \in A \bigr) + \hsp \bigr\} + \ex^{- \xx}.
\label{PcPgahGAh}
\end{EQ}
\end{theorem}

\begin{proof}
With \( \BBh = \hat{\DP} \DPc^{-1} \), \( \etav = \DPc (\vthetav - \thetavb) \), 
and \( \betav = \DPc (\hat{\thetav} - \thetavb) \)
\begin{EQA}
	\P \bigl( \hat{\DP} (\vthetav - \hat{\thetav}) \in A \cond \Yv \bigr) 
	&=&
	\P\bigl( \BBh (\etav - \betav) \in A \cond \Yv \bigr) 
	= 
	\P\bigl( \BBh (\gammav - \betav) \in A \cond \Yv \bigr) .
\label{hDPvtts}
\end{EQA}
Now the result follows from Theorem~\ref{TBvM}
and Lemma~\ref{KullbTVd} below.
\end{proof}

Therefore, the use of the estimates \( \hat{\thetav} \) and 
\( \hat{\DP}^{2} \) in place of \( \thetavb \) and \( \DPc^{2} \) 
does not significantly affect the posterior mass of any set \( A \) provided that 
the quantities \( \| \DPc (\hat{\thetav} - \thetavb) \| \) and 
\( \tr \bigl( \DPc^{-1} \hat{\DP}^{2} \DPc^{-1} - \Id_{\dimp} \bigr)^{2} \)
are small.

Theorem~\ref{TBvM} justifies the use of the posterior mean \( \vthetavb \) 
in place of \( \thetavb \).
The next important question is whether the posterior covariance 
\( \Covpost^{2} \) is a reasonable estimate of \( \DPc^{-2} \).
Unfortunately, the result \eqref{CovpostBvM} only implies  
\begin{EQA}
	\tr \bigl( \DPc^{-1} \Covpost^{2} \DPc^{-1} - \Id_{\dimp} \bigr)^{2} 
	& \leq &
	\CONST \, \dimp \, \spreads^{2}(\xx)  .
\label{trDPCovPDP}
\end{EQA}
This yields that the use of credible sets in the form
\begin{EQA}[c]
    \CR(A)
    =
    \bigl\{ \thetav \colon  \Covpost^{-1} (\thetav - \vthetavb) \in A \bigr\} 
\label{credsetBab}
\end{EQA}
is only justified if \( \dimp \, \spreads^{2}(\xx) \) is small. 
If the dimension \( \dimp \) is fixed we only need \( \spreads(\xx) \) small. 
If the dimension \( \dimp \) is large, the use of posterior covariance 
requires a stronger condition ``\( \dimp \, \spreads^{2}(\xx) \) small''.  
In the regular i.i.d. case, \( \spreads(\xx) \asymp \sqrt{(\dimp + \xx)^{3}/\nsize} \),
and \( \dimp \, \spreads^{2}(\xx) \asymp (\dimp + \xx)^{4}/\nsize \).

Alternatively, one can use a plug-in estimator of the matrix \( \DPc^{2} \).
Namely, suppose that the matrix 
\( \DP^{2}(\thetav) \eqdef - \nabla^{2} \E L(\thetav) \) at the point \( \thetav \) is 
available.
This is always the case when the model assumption \( \P \in (\P_{\thetav}) \) is correct.
Then one can define 
\( \hat{\DP}^{2} = \DP^{2}(\hat{\thetav}) \), where 
\( \hat{\thetav} \) is a pilot estimator of \( \thetavs \).
Due to Theorem~\ref{TBvM}, the posterior mean \( \vthetavb \) is a natural candidate for
\( \hat{\thetav} \) leading to the credible sets of the form
\begin{EQA}[c]
	\bar{\CR}(A)
	\eqdef
	\bigl\{ \thetav \colon \DP(\vthetavb) (\thetav - \vthetavb) \in A \bigr\},
	\qquad
	A \subset \R^{\dimp} .
\label{barCRBvM}
\end{EQA}
The only condition to check is that 
\( \tr \bigl( \DPc^{-1} \DP^{2}(\thetav) \DPc^{-1} - \Id_{\dimp} \bigr)^{2} \)
is small for all \( \thetav \) from the set \( \Thetas \) on which the estimator 
\( \vthetavb \) concentrates with a dominating probability.
The condition \( (\LL_{0}) \) implies 
\( \| \DPc^{-1} \DP^{2}(\thetav) \DPc^{-1} - \Id_{\dimp} \| \leq \rddelta(\rups) \)
for \( \thetav \in \Thetas(\rups) \) that implies
\begin{EQA}
	\tr \bigl( \DPc^{-1} \DP^{2}(\thetav) \DPc^{-1} - \Id_{\dimp} \bigr)^{2}
	& \leq &
	\dimp \, \rddelta^{2}(\rups) ,
	\quad
	\thetav \in \Thetas(\rups) .
\label{trDPcDP2DPc}
\end{EQA}
The condition \eqref{rrrupsboundco} forces \( \rups^{2} > \dimp \), 
so \( \dimp \, \rddelta^{2}(\rups) < \rups^{2} \rddelta^{2}(\rups) \leq \spreads^{2}(\xx) \).
One can see that the plug-in approach is well justified provided that the matrix function 
\( \DP^{2}(\thetav) \) is available,
in particular, under the true model specification.

\begin{corollary}
\label{CellCS}
Let the conditions of Theorem~\ref{TBvM} be fulfilled,
and in addition, the matrix function \( \DP^{2}(\thetav) \) be known.
Then Theorem~\ref{TBvMhat} applies 
to the credible sets of the form \eqref{barCRBvM} with 
\( \fiD = 1 + \rddelta(\rups) \) and \( \eps^{2} = \dimp \rddelta^{2}(\rups) \).
\end{corollary}

\Subsection{Extension to a continuous prior}
\label{Sregularprior}
The previous results for a non-informative prior can be extended to the case
of a general prior \( \Prior(d\thetav) \) with a density \( \prior(\thetav) \)
which is uniformly continuous on the local set \( \Thetas(\rups) \).
More precisely, let \( \prior(\thetav) \) satisfy
\begin{EQA}
    \sup_{\thetav \in \Thetas(\rups)}
        \bigl| \log \prior(\thetav) - \log \prior(\thetavs) \bigr|
    & \le &
    \alpha(\rups) ,
    \quad
    \sup_{\thetav \in \Theta \setminus \Thetas(\rups)} \frac{\prior(\thetav)}{\prior(\thetavs)}
    \le
    \CONST(\rups) ,
    \qquad
\label{priordenalp}
\end{EQA}
where \( \alpha(\rups) \) is a small constant while \( \CONST(\rups) \) is any fixed constant.
The second condition in \eqref{priordenalp} is not restrictive but it implicitly assumes 
that the true point \( \thetavs \) belongs to some level set of the prior density. 
Then the results of Theorem~\ref{TBvM} continue to apply to the posterior 
\( \P_{\prior}(\cdot \cond \Yv) \) 
for the prior \( \Prior \) with an obvious correction of the errors.
Indeed, for any local set \( A \subseteq \Thetas(\rups) \) one can apply the bounds
\begin{EQA}
	\int_{A} \exp\bigl\{ L(\thetav) \bigr\} \prior(\thetav) d\thetav
 	& \leq &
	\ex^{ \alpha(\rups) } \prior(\thetavs) 
	\int_{A} \exp\bigl\{ L(\thetav) \bigr\} d\thetav,
	\\
	\int_{A} \exp\bigl\{ L(\thetav) \bigr\} \prior(\thetav) d\thetav
 	& \geq &
	\ex^{ - \alpha(\rups) } \prior(\thetavs) 
	\int_{A} \exp\bigl\{ L(\thetav) \bigr\} d\thetav .
\label{priordenAl}
\end{EQA}
This particularly implies for each \( A \subset \Thetas(\rups) \) that
\begin{EQA}
	\P_{\prior}(A \cond \Yv)
	& \leq &
	\exp \bigl\{ 2 \alpha(\rups) \bigr\} \P(A \cond \Yv) .
\label{PprApo}
\end{EQA}
The tail probability of the complement \( \Thetas^{c}(\rups) \eqdef \Theta \setminus \Thetas(\rups) \)
of \( \Thetas(\rups) \) can be enlarged by \( \CONST(\rups) \)
relative to the uniform prior:
\begin{EQA}
	\int_{\Thetas^{c}(\rups)} \exp\bigl\{ L(\thetav) \bigr\} \prior(\thetav) d\thetav
 	& \leq &
	\CONST(\rups) \, \prior(\thetavs) 
	\int_{\Thetas^{c}(\rups)} \exp\bigl\{ L(\thetav) \bigr\} d\thetav,
\label{priordenTs}
\end{EQA}
hence 
\begin{EQA}
	\P_{\prior}(\Thetas^{c}(\rups) \cond \Yv)
	& \leq &
	\CONST(\rups) \P(\Thetas^{c}(\rups) \cond \Yv).
\label{PprTsc}
\end{EQA}
In particular, if the tail of the non-informative posterior satisfies 
\( \P(\Thetas^{c}(\rups) \cond \Yv) \leq \ex^{-\xx} \), then 
\( \P_{\prior}(\Thetas^{c}(\rups) \cond \Yv) \leq \CONST(\rups) \ex^{-\xx} \).

As an example, consider the case of a Gaussian prior
\( \Prior = \ND(0,\GP^{-2}) \) with the density
\( \prior(\thetav) \propto \exp \bigl\{ - \|  \GP \thetav \|^{2}/2 \bigr\} \).
The non-informative prior can be viewed as a limiting case of a Gaussian prior as
\( \GP \to 0 \).
We are interested in quantifying this relation.
How small should \( \GP \) be to ensure the BvM result?
The answer is given by the next theorem.

\begin{theorem}
\label{TGaussprior}
Suppose the conditions of Theorem~\ref{TBvM}.
Let also \( \Prior = \ND(0,\GP^{-2}) \) be a Gaussian prior measure on
\( \R^{\dimp} \) such that
\begin{EQA}[c]
    \| \DPc^{-1} \GP^{2} \DPc^{-1} \|_{\infty} \, \leq \eps^{2} ,
\label{GPthGP2ta}
\end{EQA}
where \( \eps \) is a given constant.
Then \eqref{PprApo} and \eqref{PprTsc} hold with 
\( \CONST(\rups) \leq \exp(\| \GP \thetavs \|^{2}/2) \) and 
\( \alpha(\rups) = \max\bigl\{ \eps \, \rups \| \GP \thetavs \|, \eps^{2} \rups^{2}/2 \bigr\} \). 
\end{theorem}

Note that the conditions \eqref{priordenalp} effectively mean that the prior \( \Prior \) 
does not significantly affect the posterior distribution.
Similar conditions and results can be found in the literature for more specific models, 
see e.g. \cite{Bo2011} for the Gaussian case.

\Section{Critical dimension. Examples}
\label{Scritdim}
\label{SBvMiid}
This section discusses the issue of the critical dimension of the parameter space \( \Theta \). 
This particularly allows to consider the case of a growing parameter dimension.
It appears that the error of expansions is different for different types of problems under consideration. 
The large deviation result of Theorem~\ref{TMLE} for the MLE \( \tilde{\thetav} \) requires 
\eqref{cgmibrr}.
In particular, the value \( \rups^{2} \) should be at least \( \CONST (\dimp + \xx) \). 
This result yields consistency of the MLE if the neighborhood \( \Thetas(\rups) \) is small,
or, equivalently, \( (\dimp + \xx)^{1/2} \| \DPc^{-1} \|_{\infty} \) is small. 
The Fisher expansion \eqref{ttusmxivF} of Theorem~\ref{Tconflocro} requires that 
the error term \( \Excgr(\rups,\xx) \) from \eqref{Exceqrrrhdef} is small.
Finally, the Wilks expansion and the BvM Theorem require that the error term 
\( \spread(\rups,\xx) \asymp \spreads(\xx) \) is small.
Now we specify the results for some popular statistical models.

\Subsection{Linear and generalized linear models}
In the case of a linear model 
\( \Yv = \Psi^{\T} \thetav + \epsv \) with a given design \( \dimp \times n \) matrix \( \Psi \)
under the assumption of Gaussian noise \( \epsv \sim \ND(0,\Sigma) \), the standard calculus 
leads to the log-likelihood 
\begin{EQA}[c]
	L(\thetav)
	=
	- \frac{1}{2} (\Yv - \Psi^{\T} \thetav)^{\T} \Sigma^{-1} (\Yv - \Psi^{\T} \thetav) + R ,
\end{EQA}
where the remainder \( R \) does not depend on \( \thetav \). 
Moreover, \( L(\thetav) \) is quadratic in \( \thetav \) and its Hessian 
is constant: \( \nabla^{2} \E L(\thetav) = - \Psi \Sigma^{-1} \Psi^{\T} \).
One can summarize as follows:
with \( \E \Yv = \fv \)
\begin{EQA}
	\DPc^{2} & \eqdef &
	- \nabla^{2} \E L(\thetav) 
	= 
	\Psi \Sigma^{-1} \Psi^{\T} ,
	\\
	\tilde{\thetav} 
	&=& 
	\DPc^{-2} \Psi \Sigma^{-1} \Yv,
	\\
	\thetavs 
	&=& 
	\DPc^{-2} \Psi \Sigma^{-1} \fv,
	\\
	\xiv 
	& \eqdef &
	\DPc^{-1} \nabla L(\thetavs) = \DPc^{-1} \Psi \Sigma^{-1} (\Yv - \fv).
\end{EQA}
Moreover,
\begin{EQA}
	\DPc \bigl( \tilde{\thetav} - \thetavs \bigr)
	& \equiv &
	\xiv,
	\\
	L(\tilde{\thetav}) - L(\thetavs)
	& \equiv &\
	\| \xiv \|^{2}/2 .
\label{DPcttsLLts}
\end{EQA}
All these results are straightforward, the last one is obtained by 
the Tailor expansion of the second order around \( \tilde{\thetav} \) with the use of 
\( \nabla L(\tilde{\thetav}) = 0 \):
\begin{EQA}
	L(\tilde{\thetav}) - L(\thetavs) 
	& = &
	- \frac{1}{2} (\tilde{\thetav} - \thetavs)^{\T} \nabla^{2} L(\tilde{\thetav}) 
		(\tilde{\thetav} - \thetavs)
	=
	- \frac{1}{2}\bigl\| \DPc \bigl( \tilde{\thetav} - \thetavs \bigr) \bigr\|^{2}
	=
	- \frac{1}{2} \| \xiv \|^{2} .
\end{EQA}
The presented derivations mean that 
the Fisher and Wilks expansions are \emph{identities}, they apply for \emph{any sample size}
without \emph{any conditions},
and are only based on \emph{quadraticity} of the likelihood function 
\( L(\thetav) \) in \( \thetav \).
The \emph{true distribution} of \( \Yv \) can be whatever and
it is not involved at all. 
There is \emph{no dimensional restrictions}.
However, for \emph{inference}, the parametric assumption is important.
It only concerns the \emph{distribution of \( \xiv \)}.
Let \( \Var\bigl( \Yv \bigr) = \Sigmas \ne \Sigma \).
Then with \( \DPc^{2} = \Psi \Sigma^{-1} \Psi^{\T} \)
\begin{EQA}[c]
	\Var \bigl\{ \nabla L(\thetavs) \bigr\}
	=
	\Var\bigl\{ \Psi \Sigma^{-1} \Yv \bigr\}
	=
	\Psi \Sigma^{-1} \Sigmas \Sigma^{-1} \Psi^{\T}
	\eqdef
	\VPc^{2}
	\ne 
	\DPc^{2}
\end{EQA}
which leads to the famous \emph{sandwich formula}
\begin{EQA}[c]
	\Var(\xiv)
	=
	\Var\bigl\{ \DPc^{-1} \nabla L(\thetavs) \bigr\}
	=
	\DPc^{-1} \VPc^{2} \DPc^{-1} 
	\ne 
	\Id_{\dimp} .
\end{EQA}

\Subsection{Generalized linear models (GLM)}

Let \( \Yv = (Y_{1},\ldots, Y_{n})^{\T} \sim \P \) be a sample of independent r.v.s.
The parametric GLM model is given by 
\( Y_{i} \sim P_{\Psi_{i}^{\T} \thetav} \in (P_{\ups}) \),
where \( \Psi_{i} \) are given factors in \( \R^{\dimp} \),
\( \thetav \in \R^{\dimp} \) is the unknown parameter in \( \R^{\dimp} \), and
\( (P_{\ups}) \) is an exponential family with canonical parametrization yielding the log-density
\( \ell(y,\ups) = y \ups - d(\ups) \) for a convex function \( d(\ups) \).
The MLE \( \tilde{\thetav} \) and the target \( \thetavs \) for this GLM read as
\begin{EQA}
	\tilde{\thetav}
	&=&
	\argmax_{\thetav} L(\thetav)
	=
	\argmax_{\thetav} \sum_{i=1}^{n} 
		\bigl\{ Y_{i} \Psi_{i}^{\T} \thetav - d(\Psi_{i}^{\T} \thetav) \bigr\} ,
	\\
	\thetavs
	&=&
	\argmax_{\thetav} \E L(\thetav)
	=
	\argmax_{\thetav} \sum_{i=1}^{n} 
		\bigl\{ f_{i} \Psi_{i}^{\T} \thetav - d(\Psi_{i}^{\T} \thetav) \bigr\}
\end{EQA}
with \( f_{i} = \E Y_{i} \).
An important feature of a GLM is that the stochastic component \( \zeta(\thetav) \) of 
\( L(\thetav) \) is \emph{linear in} \( \thetav \): with \( \varepsilon_{i} = Y_{i} - \E Y_{i} \)
\begin{EQA}
	\zeta(\thetav)
	&=&
	L(\thetav) - \E L(\thetav)
	=
	\Bigl( \sum_{i=1}^{n} \varepsilon_{i} \Psi_{i} \Bigr)^{\T} \thetav ,
	\quad
	\nabla \zeta(\thetav)
	=
	\sum_{i=1}^{n} \varepsilon_{i} \Psi_{i}.
\end{EQA}
%
In the contrary to the linear case, 
the Fisher information matrix \( \DPc^{2} \) depends on 
the true data distribution via the target \( \thetavs \):
\begin{EQA}[c]
	\DPc^{2}
	=
	\sum_{i} \Psi_{i} \Psi_{i}^{\T} d''(\Psi_{i}^{\T} \thetavs) .
\end{EQA}
The vector \( \xiv \) is given by
\begin{EQA}[c]
	\xiv
	=
	\DPc^{-1} \nabla \zeta(\thetavs)
	=
	\DPc^{-1} \sum_{i=1}^{n} \varepsilon_{i} \Psi_{i} \, .
\end{EQA}
Here is a list of sufficient conditions which ensure our general conditions 
from Section~\ref{Scondgllo}:
the functions \( d''(\Psi_{i}^{\T} \thetav) \) are uniformly continuous in \( \thetav \) and
\( i \leq n \),
for some matrices \( \vp_{i}^{2} \) and fixed constants \( \CONST_{0}, \lambda_{0} > 0 \)
\begin{EQA}[c]
	\E \exp\bigl\{ \lambda_{0} \vp_{i}^{-1} \varepsilon_{i} \bigr\} 
	\leq 
	\CONST_{0} ,
	\qquad 
	i=1,\ldots,\nsize,
\end{EQA}
the sample size \( \nsize \) is larger than \( \CONST \dimp \) for a prescribed constant \( \CONST \),
and the matrix \( \VPc^{2} = \sum_{i} \vp_{i}^{2} \) fulfills
\begin{EQA}[c]
	\VPc^{2}
	\leq 
	\fis^{2} \DPc^{2} .
\end{EQA}
The details can be found in \cite{SP2011}.

As already mentioned, the error of Fisher and Wilks approximations only include the term 
\( \rddelta(\rups) \) because the stochastic term is linear by definition. 
By using the higher order expansion of the function \( d(\cdot) \), 
\cite{Po1988} showed that the error \( \err(\rups,\xx) \) can be improved under the true
the parametric assumption, however, the case of a model misspecification is not include.
Under mild regularity conditions on the design \( \Psi \) and on smoothness of \( d(\cdot) \),
the consistency result applies with \( \rups^{2} = \CONST (\dimp + \xx) \), and
the bound \( \rddelta(\rups) \asymp \rups/\sqrt{n} \) yields the errors 
\( \err(\rups,\xx) \asymp \rups \rddelta(\rups) \leq \CONST (\dimp + \xx)/\sqrt{n} \) and
\( \spread(\rups,\xx) \asymp \rups^{2} \rddelta(\rups) \leq \CONST (\dimp + \xx)^{3/2}/\sqrt{n} \).

\Subsection{I.i.d. case}
\label{Chgexamples}

Now we consider the asymptotic setup in which
\( \Yv = (Y_{1},\ldots,Y_{n})^{\T} \) is an i.i.d. sample from a measure \( P \).
The parametric assumption means that \( P \in (P_{\thetav}, \,\thetav \in \Theta) \)
for a given family of marginal measures \( (P_{\thetav}) \).
We admit that the parametric family depends on \( \nsize \) and the parameter dimension 
\( \dimp = \dimn \) grows to infinity with the sample size \( \nsize \).
We also allow that the parametric assumption  can be misspecified.

Similarly to \cite{SP2011}, our general conditions can be transformed into the conditions on 
the family \( (P_{\thetav}) \) and the marginal measure \( P \).
It suffices to check that the first and second derivatives of the log-density function 
\( \ell(y,\thetav) = \log dP_{\thetav}/d\Pdom(y) \)
have exponential moments and the expectation \( \E \ell(Y_{1},\thetav) \)
is three times continuously differentiable in \( \thetav \).
See Section 5.1 in \cite{SP2011} for more details.

Also select \( \xx = \xxn \) depending on \( \nsize \) and growing slowly with \( \nsize \),
for instance, \( \xxn = \log \nsize \).
The matrix \( \DPc^{2} \) satisfies
\( \DPc^{2} = \nsize \IF_{\thetavs} \), where \( \IF_{\thetavs} \) is the Fisher information of
\( (P_{\thetav}) \) at \( \thetavs \) if the parametric assumption holds.

The bracketing bound and the large deviation result from \cite{SP2011} and from Section~\ref{Sfreqres} apply if the sample size
\( \nsize \) fulfills \( \nsize \ge \CONST \dimp_{\nsize} \) for a fixed constant
\( \CONST \).
It appears that the Fisher, Wilks, and BvM results require stronger conditions.
Indeed, in the regular i.i.d. case it holds
\( \rddelta(\rups) \asymp \rups/\sqrt{\nsize} \),
\( \qqQ^{2}(\xxn) \asymp \dimn + \xxn \),
and \( \rhor \asymp 1/\sqrt{\nsize} \).
The radius \( \rups \) should fulfill \( \rups^{2} \ge \CONST \dimp_{\nsize} \)
to ensure the large deviation result.
This yields
\begin{EQA}[c]
    \Excgr(\rups,\xxn)
    \ge
    \bigl\{ \rddelta(\rups) + 3 \nunu \qqQ(\xxn) \rhor \bigr\} \rups
    \leq 
    \CONST \dimn /  \sqrt{\nsize} .
\label{taurddimpn}
\end{EQA}
Similarly
\begin{EQA}[c]
	\spread(\rups,\xxn)
    \leq 
    \bigl\{ \rddelta(\rups) + 6 \nunu \qqQ(\xxn) \rhor \bigr\} \rups^{2} 
    \ge
    \CONST \sqrt{\dimn^{3} / \nsize} .
\label{spreadasm}
\end{EQA}
One can conclude that the consistency result is valid under 
\( \dimn/\nsize \to 0 \),
the Fisher expansion requires \( \dimn^{2} / \nsize \to 0 \),
while the Wilks and BvM results are applicable under 
\( \dimp_{\nsize}^{3}/\nsize \to 0 \) as \( \nsize \to \infty \).

\begin{theorem}
\label{TiidBvM}
Suppose the conditions of Theorem~5.4 in \cite{SP2011}.
If \( \dimn^{2} / \nsize \to 0 \), then the Fisher 
expansion \eqref{ttusmxivF} of Theorem~\ref{Tconflocro} holds with the error term 
\( \Excgr(\rups,\xxn) \to 0 \).
Let also \( \dimp_{\nsize}^{3}/\nsize \to 0 \).
Then the error term \( \spread(\rups,\xxn) \) in the Wilks expansion \eqref{WilLLe}
satisfies 
\( \spread(\rups,\xxn) \to 0 \).
The same is true for \( \spreads(\xxn) \) from Theorem~\ref{TBvM}.
\end{theorem}

Existing statistical literature addresses the issue of a growing parameter dimension 
in different set-ups. 
The classical results by \cite{Portnoy1984,Portnoy1985,Portnoy1986} 
provide some constraints on parameter dimension for consistency and asymptotic normality 
of the M-estimator for regression models. 
Our results are consistent with the conclusion of that papers. 
%
\cite{Mammen1993,Mammen1996} 
discussed the validity of bootstrap procedures in linear models with many parameters.
The obtained results are valid under \( \dimn^{3/2} / \nsize \to 0 \), however 
are limited to the testing problem for linear models. 

Another constraint in the BvM Theorem on the dimension growth \( \dimn \) can be found in \cite{Gh1999} 
for linear regression models; see the condition (2.6) 
\( \dimn^{3/2} (\log \dimn)^{1/2} \, \eta_{\nsize} \to 0 \) there, 
in  which \( \eta_{n} \) is of order \( (\dimn/\nsize)^{-1/2} \) in regular situations 
yielding a suboptimal constraint \( \nsize^{-1} \dimn^{4} \log \dimp \to 0 \).
\cite{Gh2000} obtained a version of the BvM result under the condition 
\( \nsize^{-1} \dimn^{3} (\log \dimn) \to 0 \) for a class of exponential models. 
A forthcoming paper \cite{PaSp2013} presents an example illustrating that the condition 
\( \dimp_{\nsize}^{3}/\nsize \to 0 \) cannot be dropped or relaxed.

The setup with growing parameter dimension is naturally used in sieve nonparametric estimation
when a nonparametric model is approximated by a sequence of parametric ones. 
We mention papers by 
\cite{ShWo1994, shen1997}, 
\cite{BiMa1993}, 
\cite{vdG1993,vdG2002}.
Some minimal smoothness assumptions are normally imposed on the underlying nonparametric function 
which ensure that the parameter dimension of a sieve is smaller in order than the sample size.

\Chapter{Some auxiliary results and proofs}
\label{SauxBvM}
\label{ChProofsBvM}

This section collects some auxiliary results about the behavior of the posterior measures which might be of independent interest.

\Section{Local linear approximation of the gradient of the log-likelihood}
The principle step of the proof is a bound on the local linear approximation of the gradient
\( \nabla L(\thetav) \). 
Below we study separately its stochastic and deterministic components coming from the decomposition \( L(\thetav) = \E L(\thetav) + \zeta(\thetav) \).
With \( \DPc^{2} = - \nabla^{2} \E L(\thetavs) \),
this leads to the decomposition
\begin{EQA}
\label{rderrdeftts}
	\rderr(\thetav,\thetavs)
	& \eqdef & 
	\DPc^{-1} \bigl\{ \nabla L(\thetav) - \nabla L(\thetavs) + \DPc^{2} \, (\thetav - \thetavs) \bigr\}
	\\
	&=&
	\DPc^{-1} \bigl\{ \nabla \zeta(\thetav) - \nabla \zeta(\thetavs) + 
	\nabla \E L(\thetav) - \nabla \E L(\thetavs) 
	+ \DPc^{2} \, (\thetav - \thetavs) \bigr\} .
\end{EQA}
First we check the deterministic part.
For any \( \thetav \) with \( \| \DPc (\thetav - \thetavs) \| \leq \rr \)
\begin{EQA}
	\E \bigl[ \rderr(\thetav,\thetavs) \bigr]
	& \eqdef & 
	\DPc^{-1} \bigl\{ 
		\nabla \E L(\thetav) - \nabla \E L(\thetavs) 
	\bigr\}
    + \DPc (\thetav - \thetavs) 
	\\
	&=&
	\bigl\{ \Id_{\dimp} - \DPc^{-1} \DP^{2}(\thetavd) \DPc^{-1} \bigr\} \, 
    \DPc (\thetav - \thetavs) ,
\end{EQA}
where \( \thetavd \) is a point on the line connecting \( \thetavs \) and \( \thetav \).
This implies by \( (\LL_{0}) \)
\begin{EQA}[c]
	\E \bigl[ \rderr(\thetav,\thetavs) \bigr]
	\leq 
    \| \Id_{\dimp} - \DPc^{-1} \DP^{2}(\thetavd) \DPc^{-1} \|_{\infty} \, 
    \rr
    \le 
	\rddelta(\rr) \rr .
\label{EUPupsrr}
\end{EQA}
Now we study the stochastic part. 
Consider the vector process 
\begin{EQA}
	\UP(\thetav,\thetavs)
	& \eqdef &
	\DPc^{-1} \bigl\{ 
		\nabla \zeta(\thetav) - \nabla \zeta(\thetavs)  
	\bigr\} .
\label{UPupsnm}
\end{EQA}
It is convenient to change the variable by \( \ups = \DPc (\thetav - \thetavs) \) and consider 
the vector process \( \UU(\ups) = \UP(\thetav,\thetavs) \).
It obviously holds 
\( \nabla \UU(\ups) = \DPc^{-1} \nabla^{2} \zeta(\thetav) \DPc^{-1} \). 
Moreover, for any \( \gammav_{1}, \gammav_{2} \in \R^{\dimp} \) with 
\( \| \gammav_{1} \| = \| \gammav_{2} \| = 1 \), condition \( (E\!D_{2}) \) implies
\begin{EQA}
    \log \E \exp\biggl\{ 
    	\frac{\lambda}{\rhor} \gammav_{1}^{\T} \nabla \UU(\ups) \gammav_{2} 
	\biggr\} 
    &=&
    \log \E \exp\biggl\{ 
       \frac{\lambda}{\rhor} 
       \gammav_{1}^{\T} \DPc^{-1} \nabla^{2} \zeta(\thetav) \DPc^{-1} \gammav_{2}  
    \biggr\} 
    \le 
    \frac{\nunu^{2} \lambda^{2}}{2} .
    \qquad
\label{gUUemg}
\end{EQA}
Define
\( \Upss(\rr) \eqdef \{ \ups \colon \| \ups \| \leq \rr \} \).
Then
\begin{EQA}
	\sup_{\thetav \in \Thetas(\rr)} \| \UP(\thetav,\thetavs) \|
	& = &
	\sup_{\ups \in \Upss(\rr)} \| \UU(\ups) \| .
\label{bouuvupsrupsdx1}
\end{EQA}
Theorem~\ref{Tsqnorm} 
yields
\begin{EQA}[c]
	\sup_{\ups \in \Upss(\rr)} \| \UU(\ups) \|
	\leq 
	6 \nunu \, \qqQ(\xx) \, \rhor \, \rr
\end{EQA}
on a set of a dominating probability at least \( 1 - \ex^{-\xx} \), where
the function \( \qqQ(\xx) \) is given by the following rules:
\begin{EQA}
	\qqQ(\xx)
	&=&
	\begin{cases}
		\sqrt{\QQg + 2\xx} , & 
		\text{ if } {\QQg + 2 \xx} \leq \gm^{2}, \\
		\gm^{-1} \xx + \frac{1}{2} \bigl( \gm^{-1} \QQg + \gm \bigr) , & 
		\text{ if } {\QQg + 2 \xx} > \gm^{2} .
	\end{cases}
\label{zzxxgfindef}
\end{EQA}
Here 
\( \QQg = 4 \dimp \) and 
\( \QQq = 2 \dimp^{1/2} \); see Theorem~\ref{Tsmoothpenlc} in the Appendix.

Putting together the bounds \eqref{EUPupsrr} and \eqref{bouuvupsrupsdx1} imply 
the following result.

\begin{proposition}
\label{TliapprLL}
Suppose that the matrix \( \DP^{2}(\thetav) \eqdef - \nabla^{2} \E L(\thetav) \) 
fulfills the condition \( (\LL_{0}) \) and 
let also \( (E\!D_{2}) \) be fulfilled on \( \Thetas(\rr) \) 
for any fixed \( \rr \).
Then
\begin{EQA}[c]
	\P\biggl\{ 
		\sup_{\thetav \in \Thetas(\rr)} \bigl\| \DPc^{-1} \bigl\{ 
		\nabla L(\thetav) - \nabla L(\thetavs) 
	\bigr\}
    + \DPc (\thetav - \thetavs) \bigr\|
    	\geq 
    	\Excgr(\rr,\xx) 
    \biggr\} 
    \leq
    \ex^{-\xx},
\label{supupsUPdxx}
\end{EQA}
where 
\begin{EQA}[c]
    \Excgr(\rr,\xx)
    \eqdef
    \bigl\{ \rddelta(\rr) + 6 \nunu \, \qqQ(\xx) \, \rhor \bigr\} \rr .
\label{Exceqrrrh}
\end{EQA}  
\end{proposition}

The result of Proposition~\ref{TliapprLL} can be extended to the differences
\( \UP(\thetav,\thetavd) = \DPc^{-1} \bigl\{ \nabla \zeta(\thetav) - \nabla \zeta(\thetavd) \bigr\} \):
on a set of probability at least \( 1 - \ex^{-\xx} \), it holds for any 
\( \thetav,\thetavd \in \Thetas(\rr) \) and 
\( \rderr(\thetav,\thetavd) 
= \DPc^{-1} \bigl\{ \nabla L(\thetav) - \nabla L(\thetavd) \bigr\}
+ \DPc \, (\thetav - \thetavd) \)
\begin{EQA}
	\E \bigl[ \rderr(\thetav,\thetavd) \bigr]
	& \leq &
	\rddelta(\rr) \, \| \DPc (\thetav - \thetavd) \|
	\leq
	2 \rr \, \rddelta(\rr),
	\\
	\| \UP(\thetav,\thetavd) \|
	& \leq & 
	\| \UP(\thetav,\thetavs) \| + \| \UP(\thetavd,\thetavs) \| 
	\leq
	2 \rr \, \rdomega(\xx) ,
	\\
	\bigl\| \rderr(\thetav,\thetavd) \bigr\|
	& \leq &
    2 \, \Excgr(\rr,\xx) .
\label{supupsUPd}
\end{EQA}

\subsection{Local quadratic approximation of the log-likelihood}
As the next step, 
%
%
we derive a uniform deviation bound on the error of a quadratic approximation 
\( \La(\thetav,\thetavd) = (\thetav - \thetavd)^{\T} \nabla L(\thetavd) 
    - \| \DPc (\thetav - \thetavd) \|^{2}/2 \) of 
\( L(\thetav,\thetavd) \):
\begin{EQA}
	\alp(\thetav,\thetavd)
	& \eqdef &
	L(\thetav) - L(\thetavd) 
    - (\thetav - \thetavd)^{\T} \nabla L(\thetavd) 
    + \frac{1}{2} \| \DPc (\thetav - \thetavd) \|^{2} 
    \\
    &=&
    L(\thetav,\thetavd) - \La(\thetav,\thetavd)
\label{zetatuuv}
\end{EQA}
in all \( \thetav, \thetavd \in \Thetas \), where \( \Thetas \) is some vicinity 
of a fixed point \( \thetavs \).
With \( \thetavd \) fixed, the gradient \( \nabla \alp(\thetav,\thetavd) 
\eqdef \frac{d}{d\thetav} \alp(\thetav, \thetavd) \) fulfills
\begin{EQA}[c]
	\nabla \alp(\thetav,\thetavd)
	=
	\nabla L(\thetav) - \nabla L(\thetavd) 
	+ \DPc^{2} (\thetav - \thetavd)
    =
    \DPc \, \rderr(\thetav,\thetavd) ;
\end{EQA}
cf. \eqref{rderrdeftts}.
This implies 
\begin{EQA}[c]
    \alp(\thetav,\thetavd) 
    = 
    (\thetav - \thetavd)^{\T} \nabla \alp(\thetavc,\thetavd)  ,
\label{alptts}
\end{EQA}    
where \( \thetavc \) is a point on the line connecting \( \thetav \) and \( \thetavd \).
Further, 
\begin{EQA}[c]
	\bigl| \alp(\thetav,\thetavd) \bigr|
    = 
    \bigl| (\thetav - \thetavd)^{\T} \DPc \DPc^{-1} \nabla \alp(\thetavc,\thetavd) \bigr|
    \le
    \| \DPc (\thetav - \thetavd) \| 
    \sup_{\thetavc \in \Thetas(\rr)} \bigl| \rderr(\thetavc,\thetavd) \bigr| \, .
\end{EQA}
and one can apply \eqref{supupsUPd}.
This yields the following result.

\begin{proposition}
\label{Tqapprbr}
Suppose \( (\LL_{0}) \), \( (E\!D_{0}) \), and \( (E\!D_{2}) \).
For each \( \rr \), it holds on a random set 
\( \Omega(\rr,\xx) \) of a dominating probability at least \( 1 - \ex^{-\xx} \), 
it holds with any \( \thetav, \thetavd \in \Thetas(\rr) \)
\begin{EQA}[rclcrcl]
        \frac{\bigl| \alp(\thetav,\thetavs) \bigr|}{\| \DPc (\thetav - \thetavs) \|}
    & \le &
    \Excgr(\rr,\xx) ,
    & \quad &
    \bigl| \alp(\thetav,\thetavs) \bigr|
    & \le &
    \rr \, \Excgr(\rr,\xx) ,
\label{supalp12s}
    \\
        \frac{\bigl| \alp(\thetavs,\thetav) \bigr|}{\| \DPc (\thetav - \thetavs) \|}
    & \le &
    2 \Excgr(\rr,\xx) ,
    & \quad &
    \bigl| \alp(\thetavs,\thetav) \bigr|
    & \le &
    2 \rr \, \Excgr(\rr,\xx) ,
\label{supalp12st}
    \\
        \frac{\bigl| \alp(\thetav,\thetavd) \bigr|}{\| \DPc (\thetav - \thetavd) \|}
    & \le &
    2 \Excgr(\rr,\xx) ,
    & \quad &
	\bigl| \alp(\thetav,\thetavd) \bigr|
    & \le &
    4 \rr \, \Excgr(\rr,\xx) ,
\label{supalp12}
\end{EQA}    
where \( \Excgr(\rr,\xx) \) is from \eqref{Exceqrrrh}.
\end{proposition}


\Section{Proof of Theorem~\ref{TMLE}}
By definition \( \sup_{\thetav \in \Thetas(\rups)} L(\thetav,\thetavs) \geq 0 \).
So, it suffices to check that 
\( L(\thetav,\thetavs) < 0 \) for all \( \thetav \in \Theta \setminus \Thetas(\rups) \).
The proof is based on the following bound:
for each \( \rr \)
\begin{EQA}
	\P \biggl( \sup_{\thetav \in \Thetas(\rr)}
		\bigl| 
			\zeta(\thetav,\thetavs)   
			- (\thetav - \thetavs)^{\T} \nabla \zeta(\thetavs)
		\bigr|
		\geq 
		3 \nunu \, \qqQ(\xx) \, \rhor \, \rr
	\biggr)
	& \leq &
	\ex^{-\xx} .
\label{Psuptsrrrhor}
\end{EQA}
This bound is a special case of the general result from Theorem~\ref{Tsmoothpenlc} below.
It implies by Theorem~\ref{TsuprUP} with \( \rho = 1/2 \) on a set \( \Omega(\xx) \) of 
probability at least \( 1 - \ex^{-\xx} \) that for all \( \rr \geq \rups \)
and all \( \thetav \) with \( \| \DPc (\thetav - \thetavs) \| \leq \rr \)
\begin{EQA}
	\bigl| 
			\zeta(\thetav,\thetavs)   
			- (\thetav - \thetavs)^{\T} \nabla \zeta(\thetavs)
	\bigr|
	& \leq &
	\rdomega(\rr,\xx) \, \rr,
\label{zetattsna}
\end{EQA}
where 
\begin{EQA}
	\rdomega(\rr,\xx)
	&=&
	6 \nunu \, \qqQ\bigl( \xx + \log(2\rr/\rups) \bigr) \, \rhor \, .
\label{rdomdefGP}
\end{EQA}
The use of 
\( \nabla \E L(\thetavs) = 0 \) yields
\begin{EQA}
\label{LGPLaGPrr}
	\sup_{\thetav \in \Thetas(\rr)}
		\bigl| 
			L(\thetav,\thetavs) - \E L(\thetav,\thetavs)  
			- (\thetav - \thetavs)^{\T} \nabla L(\thetavs)
		\bigr|
	& \leq &
	\rdomega(\rr,\xx) \, \rr .
\end{EQA}
By Proposition~\ref{LLbrevelocro}, the vector \( \xiv = \DPc^{-1} \nabla \zeta(\thetavs) \) 
fulfills \( \P\bigl( \| \xiv \| \geq \qq(\BB,\xx) \bigr) \leq 2 \ex^{-\xx} \).
We ignore here the negligible term of order \( \ex^{-\xxc} \).
The condition \( \| \xiv \| \leq \qq(\BB,\xx) \) implies
for \( \rr \geq \rups \) 
\begin{EQA}
	&& \nquad
	\sup_{\thetav \in \Thetas(\rr)} 
		\bigl| (\thetav - \thetavs)^{\T} \nabla L(\thetavs) \bigr|
	\\
	& \leq &
	\sup_{\thetav \in \Thetas(\rr)} 
		\| \DPc (\thetav - \thetavs) \| \times \| \DPc^{-1} \nabla L(\thetavs) \|
	= 	
	\rr \| \xiv \| 
	\leq 
	\qq(\BB,\xx) \, \rr.
\label{supDLGPsGP}
\end{EQA}
Condition \( (\LL) \) implies 
\( - 2 \E L(\thetav,\thetavs) \geq \rr^{2} \gmi(\rr) \) for each \( \thetav \) with 
\( \| \DPc (\thetav - \thetavs) \| = \rr \).
We conclude that the condition  
\begin{EQA}
	\rr \gmi(\rr)
	& \geq &
	2 \qq(\BB,\xx) + 2 \rdomega(\rr,\xx),
	\quad \rr > \rups ,
\label{rr2gmirr}
\end{EQA}
ensure \( L(\thetav,\thetavs) < 0 \) for all \( \thetav \not\in \Thetas(\rups) \) with a dominating probability.

\Section{Proof of Theorem~\ref{Tconflocro}}

Let \( \rups \) be selected to ensure that 
\( \P\bigl\{ \tilde{\thetav} \not\in \Thetas(\rups) \bigr\} \leq \ex^{-\xx} \).
Furthermore, the definition of \( \tilde{\thetav} \) yields 
\( \nabla L(\tilde{\thetav}) = 0 \) and 
\begin{EQA}[c]
	\rderr(\tilde{\thetav},\thetavs)
	=
	- \DPc^{-1} \nabla L(\thetavs) + \DPc (\tilde{\thetav} - \thetavs) .
\end{EQA}
Now Proposition~\ref{TliapprLL} implies on a set of a dominating probability
\begin{EQA}[c]
	\| \DPc (\tilde{\thetav} - \thetavs) - \xiv \|
	\leq
	\Excgr(\rr,\xx) 
\label{supupsUPdxxt}
\end{EQA}
and the assertion follows.

\Section{Proof of Theorem~\ref{TWilks2r}}
The result of Proposition~\ref{Tqapprbr} for the special case with 
\( \thetav = \thetavs \) and \( \thetavd = \tilde{\thetav} \) yields in view of 
\( \nabla L(\tilde{\thetav}) = 0 \) for \( \rr = \rups \) under the condition 
\( \tilde{\thetav} \in \Thetas(\rups) \)
\begin{EQA}[c]
    \Bigl| 
        L(\tilde{\thetav},\thetavs) 
        - \| \DPc (\tilde{\thetav} - \thetavs) \|^{2} / 2
    \Bigr|
    =
    \bigl| \alp(\thetavs,\tilde{\thetav}) \bigr|
    \le 
    2 \rups \, \Excgr(\rups,\xx) .
\label{LLttsxi}
\end{EQA}    
Further, on a set of a dominating probability, it holds 
\( \| \xiv \| \leq \qq(\BB,\xx) \);
see Proposition~\ref{LLbrevelocro}.
Now it follows from \eqref{supupsUPdxxt} that
\begin{EQA}
    && \nquad
    \bigl| \| \DPc (\tilde{\thetav} - \thetavs) \|^{2} - \| \xiv \|^{2} \bigr|
    \\
    & \le &
    2 \, \| \xiv \| \cdot \| \DPc (\tilde{\thetav} - \thetavs) - \xiv \|
    + \| \DPc (\tilde{\thetav} - \thetavs) - \xiv \|^{2}
    \\
    & \le &
    2 \, \qq(\BB,\xx) \, \Excgr(\rups,\xx) + \Excgr^{2}(\rups,\xx).
\label{ttus2xi212}
\end{EQA}    
Together with \eqref{LLttsxi}, this yields 
\begin{EQA}[c]
    \bigl| L(\tilde{\thetav},\thetavs) - \| \xiv \|^{2}/2 \bigr|
    \le 
    \bigl\{ 2 \rups + \qq(\BB,\xx) \bigr\}  \, \Excgr(\rups,\xx) 
    + \Excgr^{2}(\rups,\xx) /2.
\label{LLttsxi2}
\end{EQA}    
The error term can be improved if the squared root of the excess is 
considered.
Indeed, if \( \tilde{\thetav} \in \Thetas(\rups) \)
\begin{EQA}
    && \nquad 
    \Bigl| 
    	\bigl\{ 2L(\tilde{\thetav},\thetavs) \bigr\}^{1/2} 
		- \| \DPc (\tilde{\thetav} - \thetavs) \| 
	\Bigr|
    \leq
    \frac{\bigl| 2L(\tilde{\thetav},\thetavs) 
            - \| \DPc (\tilde{\thetav} - \thetavs) \|^{2} \bigr|}
         {\| \DPc (\tilde{\thetav} - \thetavs) \|}
    \\
    & \le & 
    \frac{2 \bigl| \alp(\tilde{\thetav},\thetavs) \bigr|}
         {\| \DPc (\tilde{\thetav} - \thetavs) \|}
    \le 
    \sup_{\thetav \in \Thetas(\rups)} 
    \frac{2 \bigl| \alp(\thetav,\thetavs) \bigr|}{\| \DPc (\thetav - \thetavs) \|}
    \le 
    2 \, \Excgr(\rups,\xx).
\label{sqLLttustu}
\end{EQA}    
The Fisher expansion \eqref{supupsUPdxxt} allows to replace here the norm 
of the standardized error \( \DPc (\tilde{\thetav} - \thetavs) \) with 
the norm of the normalized score \( \xiv \).
This completes the proof of Theorem~\ref{TWilks2r}.

\Section{Proof of Theorem~\ref{TBvM}}
The whole proof is split into few important steps.
Everywhere \( \gammav \sim \ND(0,\Id_{\dimp}) \) means a standard normal vector in 
\( \R^{\dimp} \).
For a random variable \( \eta \), denote
\begin{EQA}[c]
    \Ec \eta
    \eqdef
    \E \bigl[ \eta \cond \Yv \bigr] .
\label{Ecdefy}
\end{EQA}
Below we apply for each \( \rr \)
the bracketing bound of Proposition~\ref{Tqapprbr}: 
a random set \( \Omega(\rr,\xx) \) of probability at least \( 1 - \ex^{-\xx} \)
\begin{EQA}
	\sup_{\thetav \in \Thetas(\rr)} \bigl| L(\thetav,\thetavs) - \La(\thetav,\thetavs) \bigr| 
	& \leq &
	\rr \Excgr(\rr,\xx) 
	\quad
	\text{on }
	\Omega(\rr,\xx) ,
\label{suprrspread}
\end{EQA}
where \( \spread(\rr,\xx)\) is given by \eqref{spreadrxxdef}.
The bound from Proposition~\ref{LLbrevelocro} implies
\begin{EQA}
	\| \xiv \|	
	& \leq &
	\qq(\BB,\xx) 
	\quad
	\text{on }
	\Omega(\BB,\xx) ,
\label{xivBBxx}
\end{EQA}
on a set \( \Omega(\BB,\xx) \) of a probability at least 
\( 1 - 2 \ex^{-\xx} \).
Obviously, the probability of the overlap \( \Omega(\rr,\xx) \cap \Omega(\BB,\xx) \) is at least
\( 1 - 3 \ex^{-\xx} \).
Finally, we assume the radius \( \rups \) to be fixed which has to ensure the concentration 
of the posterior on the local set \( \Thetas(\rups) \) similarly to concentration of the MLE 
\( \tilde{\thetav} \) shown in Theorem~\ref{TMLE}.

\Section{Local Gaussian approximation of the posterior. Upper bound}
As the first step, we study the properties of
\( \DPc \bigl( \vthetav - \thetavb \bigr) \),
where
\begin{EQA}[c]
    \thetavb
    = 
    \thetavs + \DPc^{-1} \xiv
    =
    \thetavs + \DPc^{-2} \nabla L(\thetavs)
\label{thetardb}
\end{EQA}
and \( \xiv = \DPc^{-1} \nabla L(\thetavs) \).
%
For any nonnegative function \( \ff \), it holds by \eqref{suprrspread}
\begin{EQA}
    \int_{\Thetas(\rups)} \ex^{ L(\thetav,\thetavs) } \,
        f\bigl( \DPc (\thetav - \thetavb) \bigr) \, d\thetav
    & \le &
    \ex^{\spreads(\xx)}
    \int_{\Thetas(\rups)} \ex^{ \La(\thetav,\thetavs) } \,
        f\bigl( \DPc (\thetav - \thetavb) \bigr) \, d\thetav .
\label{intTheta0c}
\end{EQA}
Similarly, 
\begin{EQA}
    \int_{\Thetas(\rups)} \ex^{ L(\thetav,\thetavs) }
            f\bigl( \DPc (\thetav - \thetavb) \bigr) \, d\thetav
    & \ge &
    \ex^{- \spreads(\xx)}
    \int_{\Thetas(\rups)} \ex^{ \La(\thetav,\thetavs) }
            f\bigl( \DPc (\thetav - \thetavb) \bigr) \, d\thetav .
\label{intTheta0d}
\end{EQA}
The main benefit of these bounds is that 
\( \La(\thetav,\thetavs) \) is quadratic in \( \thetav \).
This enables to explicitly evaluate the posterior and to show that the posterior
measure is nearly Gaussian.
%
In what follows
\( \gammav \) is a standard normal vector in \( \R^{\dimp} \) independent of \( \Yv \).
Define also
\begin{EQA}
	\Indru
	& \eqdef &
	\Ind\bigl\{ \vthetav \in \Thetas(\rups) \bigr\} 
\label{Indrupsdef}
\end{EQA}
and introduce the function \( \qq(\dimp,\xx) \) which describes the quantiles of the norm 
\( \| \gammav \| \)
of a standard normal vector \( \gammav \in \R^{\dimp} \); see \eqref{zz0xxG}:
\begin{EQA}
	\qq^{2}(\dimp,\xx) 
	&=& 
	\dimp + \sqrt{6.6 \dimp \xx} \vee (6.6 \xx) .
\label{qq2pxxdef}
\end{EQA}

\begin{proposition}
\label{Tprobloc}
Suppose 
\eqref{suprrspread} for \( \rr = \rups \).
Then for any nonnegative function \( \ff(\cdot) \) on \( \R^{\dimp} \),
it holds on \( \Omega(\rups,\xx) \)
\begin{EQA}
    \Ec \bigl[ \ff\bigl( \DPc (\vthetav - \thetavb) \bigr) \Indru \bigr]
    & \le &
    \exp\bigl\{ \spreads^{+}(\xx) \bigr\} \, \E \ff(\gammav) ,
\label{PYvApi}
\end{EQA}
where 
\begin{EQA}
    \spreads^{+}(\xx)
    &=&
    2 \spreads(\xx) + \nub(\rups) ,
\label{spreadprxxdef}
	\\
\label{nubmrusp}
    \nub(\rups)
    & \eqdef &
    - \log \Pc\bigl(
        \bigl\| \gammav + \xiv \bigr\| \le \rups
    \bigr) .
\end{EQA}
If \( \rups \ge \qq(\BB,\xx) + \qq(\dimp,\xx) \), then 
on \( \Omega(\BB,\xx) \), it holds \( \nub(\rups) \leq 2 \ex^{-\xx} \) and
\( \spreads^{+}(\xx) \le 2 \spreads(\xx) + 2 \ex^{-\xx} \).
\end{proposition}

\begin{proof}
We use that
\( \La(\thetav,\thetavs)
= \xiv^{\T} \DPc (\thetav - \thetavs) - \| \DPc (\thetav - \thetavs) \|^{2} / 2 \)
is proportional to the density of a Gaussian distribution.
More precisely, define
\begin{EQA}
\label{Cdeltbd}
    \Ccb(\xiv)
    & \eqdef &
    - \| \xiv \|^{2}/2 + \log (\det \DPc) - \dimp \log(\sqrt{2\pi}) .
\end{EQA}
Then
\begin{EQA}
    \Ccb(\xiv) + \La(\thetav,\thetavs)
    &=&
    - \| \DPc (\thetav - \thetavb) \|^{2}/2
    + \log (\det \DPc) - \dimp \log(\sqrt{2\pi})
    \qquad
\label{CxivDPc}
\end{EQA}
is (conditionally on \( \Yv \)) the log-density of the normal law
with the mean \( \thetavb = \DPc^{-1} \xiv + \thetavs \) and the covariance matrix
\( \DPc^{-2} \).
Change of variables \( \uv = \DPc (\thetav - \thetavb) \) implies by
\eqref{CxivDPc} for any nonnegative function \( \ff \) that
\begin{EQA}
    && \nquad
    \int_{\Thetas(\rups)} \exp \bigl\{ L(\thetav,\thetavs) + \Ccb(\xiv) \bigr\} \,
        f\bigl( \DPc (\thetav - \thetavb) \bigr) \, d\thetav
    \\
    & \le &
    \ex^{\spreads(\xx)}
    \int \exp \bigl\{ \La(\thetav,\thetavs) + \Ccb(\xiv) \bigr\} \,
        f\bigl( \DPc (\thetav - \thetavb) \bigr) \, d\thetav
    \\
    &=&
    \ex^{\spreads(\xx)} \int \phi(\uv) \, \ff(\uv) \, d\uv
    = \ex^{\spreads(\xx)} \,\, \E \ff(\gammav) .
\label{intALAc}
\label{PhimPhib}
\end{EQA}
Similarly, for any nonnegative function \( \ff \), it follows
by change of variables \( \uv = \DPc (\thetav - \thetavb) \) and 
\( \DPc (\thetav - \thetavs) = \uv + \xiv \)
that
\begin{EQA}
    && \nquad
    \int \exp \bigl\{ L(\thetav,\thetavs) \bigr\} \,
    f\bigl( \DPc (\thetav - \thetavb) \bigr) \, 
        \Ind\bigl\{ \| \DPc (\thetav - \thetavs) \| \le \rups \bigr\} d\thetav
    \\
    & \ge &
    \exp\{ - \spreads(\xx) - \Ccm(\xiv) \} \,
    \int \phi(\uv) \ff(\uv)
    \Ind\bigl\{ \| \uv + \xiv \| \le \rups \bigr\} d \uv.
    \qquad
\label{PhimPhim}
\end{EQA}
A special case of 
\eqref{PhimPhim} with
\( \ff(\uv) \equiv 1 \) implies by definition of \( \nub(\rups) \):
\begin{EQA}
    \int_{\Thetas(\rups)} \exp\{ L(\thetav,\thetavs) \} \, d \thetav
    & \ge &
    \exp \bigl\{ - \spreads(\xx) - \Ccm(\xiv) - \nub(\rups) \bigr\} .
\label{Pxivlm}
\end{EQA}
%
Further,
\eqref{PhimPhib} and \eqref{Pxivlm} imply on \( \Omega(\rups,\xx) \)
\begin{EQA}
    \frac{\int_{\Thetas(\rups)} \exp \bigl\{ L(\thetav,\thetavs) \bigr\}
            f\bigl( \DPc (\thetav - \thetavb) \bigr) \, d\thetav}
         {\int \exp \bigl\{ L(\thetav,\thetavs) \bigr\} \, d\thetav}
    & \le &
    \exp\bigl\{ 2 \spreads(\xx) + \nub(\rups)
    \bigr\} \,\, \E \ff(\gammav) 
\label{PYvApib}
\end{EQA}
and \eqref{PYvApi} follows.
As \( \| \xiv \| \leq \qq(\BB,\xx) \) on \( \Omega(\BB,\xx) \) and 
\( \rups \ge \qq(\BB,\xx) + \qq(\dimp,\xx) \),
this and Lemma~\ref{LPGxx} imply for \( \gammav \sim \ND(0,\Id_{\dimp}) \)
\begin{EQA}
    \nub(\rups)
    &=&
    - \log \Pc\bigl(
        \bigl\| \gammav + \xiv \bigr\| \leq \rups
    \bigr)
    \le 
    - \log \P\bigl( \| \gammav \| \leq \qq(\dimp,\xx) \bigr)
    \le
    2 \ex^{-\xx} ,
\label{ombrue}
\end{EQA}
and the last assertion follows.
\end{proof}

The condition ``\( \spreads(\xx) \) is small''
allows us to ignore the exp-factor in \eqref{PYvApi} and this result
yields an upper bound \( \E \ff(\gammav) \)
for the posterior expectation of
\( \ff\bigl(  \DPc (\vthetav - \thetavb) \bigr) \) conditioned on \( \Yv \)
and on \( \vthetav \in \Thetas(\rups) \).

%
The next result considers some special cases of \eqref{PYvApi} with
\( \ff(\uv) = \exp(\lambdav^{\T} \uv) \)
and 
\( \ff(\uv) = \Ind(\uv \in A) \) for a measurable subset \( A \subset \R^{\dimp} \).

\begin{corollary}
\label{TLaplaceup}
\label{CprobUBvM}
Suppose 
\eqref{suprrspread} for \( \rr = \rups \).
For any \( \lambdav \in \R^{\dimp} \), it holds on \( \Omega(\rups,\xx) \)
\begin{EQA}
\label{LapYvlamb}
    \log \Ec \bigl[ \exp\bigl\{ \lambdav^{\T} \DPc (\vthetav - \thetavb) \bigr\} \Indru \bigr]
    & \le &
    \| \lambdav \|^{2}/2 + \spreads^{+}(\xx) .
\end{EQA}
For any measurable set \( A \), it holds on \( \Omega(\rups,\xx) \)
\begin{EQA}
\label{PcPgamG}
	\Pc\bigl( \DPc (\vthetav - \thetavb) \in A \bigr)
	\eqdef
	\P\bigl( \DPc (\vthetav - \thetavb) \in A \cond \Yv \bigr)
	& \leq &
	\exp \bigl\{ \spreads^{+}(\xx) \bigr\} 
	\P\bigl( \gammav \in A \bigr) .
	\qquad
\label{PcPgamd}
\end{EQA}
\end{corollary}

In the next result we describe the local concentration properties of the posterior.
Namely, 
the centered and scaled posterior vector  
\( \DPc (\vthetav - \thetavb) \) concentrates on a coronary set 
\( \bigl\{ \uv \colon \qq_{1}(\dimp,\xx) \leq \| \uv \| \leq \qq(\dimp,\xx) \bigr\} \)
with \( \Pc \)-probability of order \( 1 - 2 \ex^{- \xx} \); see \eqref{Pgammxx}.

\begin{corollary}
\label{CetaconcBvM}
For \( \xx \geq 0 \), with \( \qq(\dimp,\xx) \) and \( \qq_{1}(\dimp,\xx) \) from \eqref{zz0xxG}
\begin{EQA}[ccl]
\label{PDPbvttBvM}
	\Pc\bigl\{ 
		\| \DPc (\vthetav - \thetavb) \| \geq \qq(\dimp,\xx)
	\bigr\}
	& \leq &
	\exp\bigl\{ - \xx + \spreads^{+}(\xx) \bigr\} ,
	\\
	\Pc\bigl( 
		\| \DPc (\vthetav - \thetavb) \| 
		\leq 
		\qq_{1}(\dimp,\xx)
	\bigr)
	& \leq &
	\exp\bigl\{ - \xx + \spreads^{+}(\xx) \bigr\} .
\label{PDPbvttBvMd}
\end{EQA}
\end{corollary}

\begin{proof}
This result is a combination of the bound from Corollary~\ref{CprobUBvM} and 
the bounds for the standard normal distributions from Lemma~\ref{LPGxx}.
\end{proof}

\Subsection{Tail posterior probability and contraction}

The next important step in our analysis is to check that \( \vthetav \) concentrates
in a small vicinity \( \Thetas = \Thetas(\rups) \) of the point \( \thetavs \)
with a properly selected \( \rups \).
This will be described by using the random
quantity
\begin{EQA}
    \rho(\rups)
    & \eqdef &
    \frac{\int_{\Theta \setminus \Thetas} \exp \bigl\{ L(\thetav) \bigr\} d \thetav}
         {\int_{\Thetas} \exp \bigl\{ L(\thetav) \bigr\} d \thetav}
    =
    \frac{\int_{\Theta \setminus \Thetas} \exp \bigl\{ L(\thetav,\thetavs) \bigr\} d \thetav}
         {\int_{\Thetas} \exp \bigl\{ L(\thetav,\thetavs) \bigr\} d \thetav}
    \, .
\label{rhopipopr}
\end{EQA}
Obviously
\( \P\bigl\{ \vthetav \not\in \Thetas(\rups) \cond \Yv \bigr\} \le \rho(\rups) \).
Therefore, small values of \( \rho(\rups) \) indicate a small posterior
probability of the set \( \Theta \setminus \Thetas \).
The proof only uses condition \( (\LL) \) and the fact that 
there exists a random set \( \Omega(\xx) \) of probability at least \( 1 - \ex^{-\xx} \) such that
\begin{EQA}
	\bigl| \zeta(\thetav,\thetavs) - \xiv^{\T} \DPc (\thetav - \thetavs) \bigr| 
	& \leq & 
	\rr \, \rdomega(\rr,\xx)
	\quad
	\text{on } \Omega(\xx)
\label{zettsxivDrr}
\end{EQA}
for \( \rr = \| \DPc (\thetav - \thetavs) \| \) 
and  \( \rdomega(\rr,\xx) \) from \eqref{rdomdefGP}; cf. the proof of Theorem~\ref{TMLE}.

Let \( \gmi_{0} = \gmi(\rups) \) and for the sequence \( \gmi_{k} = 2^{-k} \gmi_{0} \), 
the radii \( \rr_{0} < \rr_{1} < \ldots \) 
be defined by the condition 
\( \gmi(\rr) \geq \gmi_{k} > 0 \) for \( \rr_{k} \leq \rr < \rr_{k+1} \) for all
\( k \geq 0 \) with \( \gmi(\rr) \) from \( (\LL) \).

\begin{proposition}
\label{Tpostrtail}
Suppose the conditions \( (\LL) \), \( (E\!D_{0}) \), and \( (E\!D_{2}) \).
If \( \gmi(\rr) \) from \( (\LL) \) satisfies
\begin{EQA}
	\rr^{2} \gmi^{2}(\rr)
	& \geq &
	\xx
	+
	2 \dimp 
	+ 4 \qq^{2}(\BB,\xx) 
	+ 8 \rr \, \gmi(\rr) \rdomega(\rr,\xx),
	\qquad
	\rr \geq \rups,
\label{rrboundco}
\end{EQA}
then with \( \spreads^{+}(\xx) \) from \eqref{spreadprxxdef}, it holds on a set \( \Omega_{1}(\xx) \) of probability at least \( 1 - 4 \ex^{-\xx} \)
\begin{EQA}[c]
    \rho(\rups)
    \le
    2 \exp\{ - \xx + \spreads^{+}(\xx) \} \,\,
	\qquad
\label{rhointTTspr}
\end{EQA}
\end{proposition}

\begin{remark}
Suppose that \( \gmi_{0} = \gmi(\rups) \) is close to one.
Condition \eqref{rrboundco} requires in particular that 
\( \rups^{2} > 4 \qq^{2}(\BB,\xx) + 2 \dimp + \xx \) and the value \( \rr \gmi(\rr) \) grows  
with \( \rr \).
\end{remark}

\begin{proof}
For the denominator of \( \rho(\rups) \) we apply the lower bound \eqref{Pxivlm}.
It remains to bound from above the integral over the complement of the local set 
\( \Thetas(\rups) \).
Similarly to the proof of Theorem~\ref{TMLE}, 
we use the decomposition 
\( L(\thetav,\thetavs) = \E L(\thetav,\thetavs) + \zeta(\thetav,\thetavs) \).
Condition \( (\LL) \) for the expected negative log-likelihood implies
\begin{EQA}
	- \E L(\thetav,\thetavs) 
	& \geq & 
	\bigl| \DPc (\thetav - \thetavs) \bigr|^{2} \gmi_{k}/2 
\label{ELzetattsrr}
\end{EQA}
for each \( k \geq 0 \) and any \( \thetav \in \Thetas(\rr_{k+1}) \setminus \Thetas(\rr_{k}) \).
The bound \eqref{zettsxivDrr} implies on \( \Omega(\xx) \) 
\begin{EQA}
	\bigl| \zeta(\thetav,\thetavs) - \xiv^{\T} \DPc (\thetav - \thetavs) \bigr|
	& \leq &
	\rr_{k+1} \, \rdomega(\rr_{k+1},\xx),
	\qquad
	\thetav \in \Thetas(\rr_{k+1}) \setminus \Thetas(\rr_{k}) ,
\end{EQA}
for all \( k \geq 0 \).
By change of variables \( \gammav = \DPc (\thetav - \thetavs) \), it follows 
for each \( k \)  
\begin{EQA}
	&& \nquad
	\exp\bigl\{ \Ccm(\xiv) \bigr\}
	{\int_{\Thetas(\rr_{k+1}) \setminus \Thetas(\rr_{k})} 
		\exp\bigl\{ L(\thetav,\thetavs) \bigr\} d\thetav}
	\\
	& \leq &
	\exp\Bigl\{ 
		\rr_{k+1} \, \rdomega(\rr_{k+1},\xx) - \frac{\| \xiv \|^{2}}{2} 
	\Bigr\}
	\frac{1}{(2\pi)^{\dimp/2}}
	\int_{\| \gammav \| \geq \rr_{k}} \exp\bigl\{ 
		- \frac{\gmi_{k} \| \gammav \|^{2}}{2} + \xiv^{\T} \gammav 
	\bigr\} d\gammav  \, .
\label{rhoruxx}
\end{EQA}
Next, 
\begin{EQA}
	&&
	\nquad
	\frac{1}{(2\pi)^{\dimp/2}}
	\int_{\| \gammav \| \geq \rr_{k}} \exp\Bigl(
		- \frac{\gmi_{k} \| \gammav \|^{2}}{2} + \xiv^{\T} \gammav 
	\Bigr) d\gammav 
	\\
	& \leq &
	\gmi_{k}^{- \dimp/2} \exp\Bigl( \frac{\| \xiv \|^{2}}{2 \gmi_{k}} \Bigr)
	\Pc\bigl( \| \gammav + \gmi_{k}^{-1/2} \xiv \| \geq \gmi_{k}^{1/2} \rr_{k} \bigr)
	\\
	& \leq &
	\gmi_{k}^{- \dimp/2} \exp\Bigl( 
		\frac{\| \xiv \|^{2}}{\gmi_{k}} - \frac{1}{4} \gmi_{k} \rr_{k}^{2} 
		+ \frac{\dimp}{2} 
	\Bigr) .
\label{rhoruxx2}
\end{EQA}
Here we have used the bound \eqref{gampxiz12} for a standard normal vector \( \gammav \) and 
\( \uv = \gmi_{k}^{-1/2} \xiv \in \R^{\dimp} \).
\eqref{Pxivlm} and \eqref{rhoruxx2} imply \eqref{rhointTTspr}.
Now the bound 
\( \| \xiv \| \leq \qq(\BB,\xx) \) holding with a dominating probability and 
\eqref{rrboundco} imply 
\begin{EQA}
	&& \nquad
	\sum_{k=0}^{\infty} \exp\bigl\{ \Ccm(\xiv) \bigr\}
	{\int_{\Thetas(\rr_{k+1}) \setminus \Thetas(\rr_{k})} 
		\exp\bigl\{ L(\thetav,\thetavs) \bigr\} d\thetav}
	\\
	& \leq &
	\sum_{k=0}^{\infty} \exp\Bigl( 
		\frac{\| \xiv \|^{2}}{\gmi_{k}} - \frac{1}{4} \gmi_{k} \rr_{k}^{2} 
		+ \frac{\dimp}{2} \log\bigl( e / \gmi_{k} \bigr)
		+ \rr_{k+1} \, \rdomega(\rr_{k+1},\xx)
	\Bigr) 
	\\
	& \leq &
	\sum_{k=0}^{\infty} \exp(- \xx/\gmi_{k})
	\leq 
	2 \ex^{- \xx }
\label{rhosumTTspr}
\end{EQA}
and \eqref{rhointTTspr} follows in view of \( \gmi \log(e/\gmi) \leq 1 \) for \( \gmi \leq 1 \).
\end{proof}

\begin{proposition}
\label{CTpostrtail}
Assume the conditions of Proposition~\ref{Tpostrtail}.
It holds on a set \( \Omega_{2}(\xx) \) of probability at least \( 1 - 4 \ex^{-\xx} \)
for any unit vector \( \av \in \R^{\dimp} \)
\begin{EQA}
    \rho_{2}(\rups)
    & \eqdef &
    \frac{\int_{\Theta \setminus \Thetas} 
    		| \av^{\T} \DPc (\thetav - \thetavs) |^{2} 
    		\ex^{ L(\thetav,\thetavs) } d \thetav}
         {\int_{\Thetas} \ex^{ L(\thetav,\thetavs) } d \thetav}
    \leq 
    2 \exp\bigl\{ - \xx  + \spreads^{+}(\xx) \bigr\} .
    \qquad
\label{rhomrups}
\end{EQA}
\end{proposition}

\begin{proof}
The arguments are similar to the proof of Proposition~\ref{Tpostrtail} with the use of 
\eqref{gampxiz122} in place of \eqref{gampxiz12}.
\end{proof}

\Subsection{Local Gaussian approximation of the posterior. Lower bound}

Now we present a local lower bound for the posterior probability.
The reason for separating the upper and lower bounds is that the lower bound also requires
a tail probability estimation; see \eqref{rhointTTspr} and \eqref{rhomrups}.

\begin{proposition}
\label{Tproblocm}
Suppose \eqref{suprrspread} for \( \rr = \rups \) and \eqref{rhointTTspr}.
Then for any nonnegative function \( \ff(\cdot) \) on \( \R^{\dimp} \),
it holds on \( \Omega(\xx) \)
\begin{EQA}
    \Ec \bigl\{ \ff\bigl( \DPc (\vthetav - \thetavb) \bigr) \Indru \bigr\}
    & \ge &
    \exp\bigl\{ - \spreads^{-}(\xx) \bigr\} \,
        \E \Bigl\{ \ff(\gammav) \Ind\bigl( \| \gammav + \xiv \| \le \rups \bigr) \Bigr\} ,
    \qquad
\label{PYvApilb}
\end{EQA}
where \( \spreads^{-}(\xx) = \spreads^{+}(\xx) + \rho(\rups) \).
\end{proposition}

\begin{proof}
On the set \( \Omega(\xx) \), it holds by \eqref{PhimPhib} with
\( \ff(\cdot) = 1 \): 
\begin{EQA}
    \int \exp \bigl\{ L(\thetav,\thetavs) \bigr\} \, d\thetav
    & \le &
    \int_{\Thetas} \exp \bigl\{ L(\thetav,\thetavs) \bigr\} \, d\thetav
    + \int_{\Theta \setminus \Thetas} \exp \bigl\{ L(\thetav,\thetavs) \bigr\} \, d \thetav
    \\
    & \le &
    \bigl\{ 1 + \rho(\rups) \bigr\}
    \int_{\Thetas} \exp \bigl\{ L(\thetav,\thetavs) \bigr\} \, d\thetav
    \\
    & \le &
    \bigl\{ 1 + \rho(\rups) \bigr\}
    \exp\bigl\{ \spreads(\xx) - \Ccb(\xiv) + \nub(\rups) \bigr\}
    \\
    & \le &
    \exp\bigl\{ \spreads(\xx) - \Ccb(\xiv) + \nub(\rups) + \rho(\rups) \bigr\}.
\label{intTSTTS}
\end{EQA}
This and the bound \eqref{PhimPhim} imply
\begin{EQA}
    && \nquad
    \frac{\int_{\Thetas(\rups)} \exp \bigl\{ L(\thetav,\thetavs) \bigr\}
            \ff\bigl( \DPc (\thetav - \thetavb) \bigr) \, d\thetav}
         {\int \exp \bigl\{ L(\thetav,\thetavs) \bigr\} \, d\thetav}
    \\
    & \ge &
    \frac{\exp\bigl\{ - \spreads(\xx) - \Ccm(\xiv) \bigr\}
            \int \phi(\uv) \ff(\uv) \Ind\bigl\{ \| \uv + \xiv \| \leq \rups \bigr\} d\uv}
         {\exp\bigl\{ \spreads(\xx) - \Ccb(\xiv) + \nub(\rups) + \rho(\rups) \bigr\}}
    \\
    & \ge &
    \exp\bigl\{- \spreads^{-}(\xx) \bigr\} \,
        \E \bigl[ \ff(\gammav) \Ind\bigl\{ \| \gammav + \xiv \| \leq \rups \bigr\} \bigr] .
\label{PYvApim}
\end{EQA}
This yields \eqref{PYvApilb}.
\end{proof}

Note that the bound \( \| \xiv \| \leq \qq(\BB,\xx) \)
implies for \( \rups > \qq(\BB,\xx) \)
\begin{EQA}
    \bigl\{ \uv \in \R^{\dimp} \colon \| \uv + \xiv \| \leq \rups \bigr\}
    & \supseteq &
    \bigl\{
        \uv \in \R^{\dimp} \colon \| \uv \| \le \rups - \qq(\BB,\xx)
    \bigr\} .
\label{Brdrupsuv2}
\end{EQA}

As a corollary, we state the results for the distribution and moment generating functions of
\( \DPc (\vthetav - \thetavb) \).
We assume that the \( \rups^{2} \) be selected to ensure the tail probability bound \eqref{rhomrups}.

\begin{corollary}
\label{CTproblocm}
Suppose \eqref{zettsxivDrr} and \eqref{rhointTTspr}.
Let \( \rups^{2} \) ensure \eqref{rrboundco}.
On \( \Omega(\xx) \cap \Omega(\BB,\xx) \), it holds
for any \( \lambdav \in \R^{\dimp} \) with \( \| \lambdav \|^{2} \le \dimp \)
\begin{EQA}
    \log \Ec \bigl[ \exp\bigl\{ \lambdav^{\T} \DPc (\vthetav - \thetavb) \bigr\} \Indru \bigr]
    & \ge &
    \| \lambdav \|^{2}/2 - \spreads^{-}(\xx) - 2 \ex^{-\xx} .
\label{PsiYvApilm}
\end{EQA}
Moreover, for any \( A \subset \R^{\dimp} \), it holds  on \( \Omega(\xx) \) 
\begin{EQA}
	\Pc\bigl( \DPc (\vthetav - \thetavb) \in A \bigr)
	& \geq &
	\exp \bigl\{ \spreads^{-}(\xx) \bigr\} 
	\P\bigl( \gammav \in A \bigr) - \ex^{-\xx} .
\label{PcPgaGA}
\end{EQA}
\end{corollary}

\begin{proof}
The first result follows from Proposition~\ref{Tproblocm}.
The only important additional step is an evaluation of the integral
\( \E \bigl\{ \exp (\lambdav^{\T} \gammav) \Ind( \| \gammav \| \le \rr) \bigr\} \).
The bound \eqref{exIdgamrr} yields \eqref{PsiYvApilm} in view of
\( \log(1 - \ex^{-3\xx/2}) \ge - \ex^{-\xx} \) for \( \xx \ge 1 \).
The second statement can be proved similarly to Corollary~\ref{TLaplaceup}.
\end{proof}

\Subsection{Moments of the posterior}
\label{SproofBvM}
Here we show that the first two moments of the posterior are pretty close to the moments of the standard normal law.
The results are entirely based on our obtained statements from Corollaries~\ref{TLaplaceup}  
and \ref{CTproblocm}.
Due to our previous results, it is convenient to decompose the r.v. 
\( \etav = \DPc (\vthetav - \thetavb) \)
in the form
\begin{EQA}[c]
    \etav
    =
    \etav \Indru + \,
    \etav \Ind(\vthetav \not\in \Thetas(\rups))
    =
    \etavd + \etav^{c}.
\label{vthetavTTc}
\end{EQA}
The large deviation result yields that the posterior distribution of the part
\( \etav^{c} \) is negligible provided a proper choice of \( \rups \).
Below we show that \( \etavd \) is nearly standard normal which yields the
BvM result.
%
Define also the first two moments of \( \etavd \):
\begin{EQA}[c]
    \etavb
    \eqdef
    \Ec \etavd,
    \qquad
    \Covd^{2}
    \eqdef
    \Ec \bigl\{ (\etavd - \etavb) (\etavd - \etavb)^{\T} \bigr\} .
\label{meanVarvthetav}
\end{EQA}
Similarly to the proof of Corollaries~\ref{TLaplaceup} and \ref{CTproblocm}
one derives for any unit vector \( \uv \in \R^{\dimp} \)
\begin{EQA}[rcccl]
\label{Eclametav2}
    \exp \spreads^{-}(\xx)
    & \le &
    \Ec \bigl| \uv^{\T} \etavd \bigr|^{2}
    & \le &
    \exp \spreads^{+}(\xx);
\end{EQA}
see \eqref{PYvApi}, \eqref{rhomrups}, and \eqref{PYvApilb}.
%
It suffices to show that \eqref{Eclametav2} implies
\begin{EQA}[c]
    \| \etavb \|^{2} \le 2 \spreads^{*}(\xx),
    \qquad
    \| \Covd^{2} - \Id_{\dimp} \|_{\infty}
    \le
    2 \spreads^{*}(\xx) 
\label{Eexplametb}
\end{EQA}
with \( \spreads^{*}(\xx) = \max\bigl\{ \spreads^{+}(\xx), \spreads^{-}(\xx) \bigr\} \le 1/2 \).
Note now that
\begin{EQA}[c]
    \Ec \bigl| \uv^{\T} \etavd \bigr|^{2}
    =
    \uv^{\T} \Covd^{2} \uv + |\uv^{\T} \etavb|^{2} .
\label{EcuvTetavb}
\end{EQA}
Hence
\begin{EQA}[c]
    \exp\bigl\{ - \spreads^{-}(\xx) \bigr\}
    \le
    \uv^{\T} \Covd^{2} \uv + |\uv^{\T} \etavb|^{2}
    \le
    \exp \bigl\{ \spreads^{+}(\xx) \bigr\}  .
\label{spreadpmBvM}
\end{EQA}
In a similar way with \( \uv = \etavb / \| \etavb \| \) and
\( \gammav \sim \ND(0,\Id_{\dimp}) \)
\begin{EQA}
    \uv^{\T} \Covd^{2} \uv
    &=&
    \Ec \bigl| \uv^{\T} (\etav - \etavb) \bigr|^{2}
    \\
    & \ge &
    \exp\bigl\{ - \spreads^{-}(\xx) \bigr\} \E \bigl| \uv^{\T} (\gammav - \etavb) \bigr|^{2}
    =
    \exp\bigl\{ - \spreads^{-}(\xx) \bigr\} \bigl( 1 + \| \etavb \|^{2} \bigr)
\label{Ecetavvb}
\end{EQA}
yielding
\begin{EQA}[c]
    \uv^{\T} \Covd^{2} \uv
    \ge
    \bigl( 1 + \| \etavb \|^{2} \bigr) \exp \bigl\{ \spreads^{-}(\xx) \bigr\} .
\label{uvCovd2uv}
\end{EQA}
This inequality contradicts \eqref{spreadpmBvM} if
\( \| \etavb \|^{2} > 2 \spreads^{*}(\xx) \) for \( \spreads^{*}(\xx) \le 1/2 \),
and \eqref{Eexplametb} follows.

The bound for the first moment implies 
\begin{EQA}[c]
    \bigl\| \DPc (\vthetavd - \thetavb) \bigr\|^{2}
    \le
    2 \spreads^{*}(\xx)
\label{DPmEcvtrdm}
\end{EQA}
while the second bound yields \(
    \bigl\| \DPc \Covpostd^{2} \DPc - \Id_{\dimp} \bigr\|_{\infty}
    \le
    2 \spreads^{*}(\xx) \).
This completes the proof of \eqref{CovpostBvM}.

\subsection{Proof of Theorem~\ref{TGaussprior}}
It suffices to check \eqref{priordenalp}.
First evaluate the ratio 
\( \prior(\thetav)/\prior(\thetavs) \) 
for any \( \thetav \in \Thetas(\rups) \). 
It holds
\begin{EQA}[c]
    \log \frac{\prior(\thetav)}{\prior(\thetavs)}
    =
    - \| \GP \thetav \|^{2} / 2 + \| \GP \thetavs \|^{2} / 2
    =
    - (\thetav - \thetavs)^{\T} \GP^{2} \thetavs
    - \| \GP (\thetav - \thetavs) \|^{2}/2 .
\label{logpriorgauss}
\end{EQA}
%
It follows from the definition of \( \Thetas(\rups) \) and  \eqref{GPthGP2ta} that for 
\( \thetav \in \Thetas(\rups) \)
\begin{EQA}[c]
    \| \GP (\thetav - \thetavs) \|^{2}
    =
    \| \GP \DPc^{-1} \DPc (\thetav - \thetavs) \|^{2}
    \le
    \| \DPc^{-1} \GP^{2} \DPc^{-1} \|_{\infty} \, \rups^{2} 
    \leq 
    \eps^{2} \, \rups^{2}\, .
\label{GPttsVP1}
\end{EQA}
Similarly
\begin{EQA}[c]
    \bigl| (\thetav - \thetavs)^{\T} \GP^{2} \thetavs \bigr|
    \le
    \| \GP \thetavs \| \cdot \| \GP (\thetav - \thetavs) \|
    \le
    \| \GP \thetavs \| \cdot \| \GP \DPc^{-1} \|_{\infty} \, \rups 
    \leq 
    \eps \, \rups \| \GP \thetavs \| .
\label{ttvsGPtvs}
\end{EQA}
This obviously implies
\begin{EQA}
	- \frac{1}{2}\eps^{2} \, \rups^{2}
	\leq 
	\log \frac{\prior(\thetav)}{\prior(\thetavs)}
	& \leq &
	\eps \, \rups \| \GP \thetavs \|
\label{12epsruC}
\end{EQA}
and \eqref{priordenalp} follows with
\( \alpha(\rups) = \max\bigl\{ \eps \, \rups \| \GP \thetavs \|, \eps^{2} \rups^{2}/2 \bigr\} \).


\appendix

\section{An entropy bound for the maximum of a random process}
We use one general result on the upper bound for the maximum of a centered random process 
in the form of \cite{SP2011}; see Corollary~7.2 in the supplement of that paper.

Here we discuss the special case when \( \Ups \) is an open subset in 
\( \R^{\dimp} \),
the stochastic process \( \UP(\ups) \) is absolutely continuous and its gradient
\( \nabla \UP(\ups) \eqdef d \UP(\ups) / d \ups \) 
has bounded exponential moments.

\begin{description}
\item[\( \bb{(\CS\! D)} \)]\textit{
There exist \( \gmb > 0 \), 
\( \nunu \ge 1 \), and 
a symmetric non-negative matrix \( \VVc \)
such that for any  \( \lambda \le \gmb \) 
and any unit vector \( \gammav \in \R^{\dimp} \), it holds
}
\begin{EQA}[c]
    \log \E \exp \Bigl\{
       \lambda 
       \frac{\gammav^{\T}\nabla \UP(\ups)}
            {\| \VVc \gamma \|}
    \Bigr\} 
    \le 
    \frac{\nunu^{2} \lambda^{2}}{2} .
\end{EQA}
\end{description}

%


We consider the local sets of the elliptic form 
\( \Upss(\rr) \eqdef \{ \ups: \| \VVc(\ups - \upss) \| \le \rr \} \).

\begin{theorem}[\cite{SP2011}]
\label{Tsmoothpenlc}
Let \( (\CS\! D) \) hold with some \( \gmb > 0 \), and a matrix \( \VVc \).
For any \( \xx \ge 0 \) and any \( \rr > 0 \)
\begin{EQA}[c]
\label{Upsdbounddsm}
    \P \Bigl\{ 
        \sup_{\ups \in \Upss(\rr)} \bigl| \UP(\ups,\upss) \bigr| 
    	\geq 3 \nunu \, \rr \, \qqQ(\xx)
    \Bigr\}
    \le 
    \ex^{-\xx} ,
\end{EQA}
where \( \qqQ(\xx) \) is given by the following rule:
\begin{EQA}
	\qqQ(\xx)
	&=&
	\begin{cases}
		\sqrt{\QQ + 2\xx} & 
		\text{ if } {\QQ + 2 \xx} \leq \gm^{2}, \\
		\gm^{-1} \xx + \frac{1}{2} \bigl( \gm^{-1} \QQ + \gm \bigr) & 
		\text{ if } {\QQ + 2 \xx} > \gm^{2} ,	
	\end{cases}
\label{zzxxgfin}
\end{EQA}
with \( \QQ = 4 \dimp \).
\end{theorem}

Due to the result of Theorem~\ref{Tsmoothpenlc}, the bound for the maximum of 
\( \UP(\ups,\upss) \) over \( \ups \in \B_{\rr}(\upss) \) grows linearly in \( \rr \). 
So, its applications to situations with \( \rr \gg \entrlq(\Upsd) \) are limited.
The next result shows that introducing a negative drift helps to
state a uniform in \( \rr \) local probability bound.
Namely, the bound for the process 
\( \UP(\ups,\upss) - f(\dist(\ups,\upss)) \) for some 
function \( f(\rr) \) over a ball \( \B_{\rr}(\upss) \) around 
the point \( \upss \) does not depend on \( \rr \).
Here the generic chaining arguments are accomplished with the slicing technique.
The idea is for a given \( \rrb > 1 \) to split the ball \( \B_{\rrb}(\upss) \) into the slices 
\( \B_{\rr_{k}}(\upss) \setminus \B_{\rr_{k-1}}(\upss) \) 
and to apply Theorem~\ref{Tsmoothpenlc} to each slice separately.

\begin{theorem}
\label{TsuprUP}
Let \( \rrb \) be such that \( (\CS{D}) \) holds on \( \B_{\rrb}(\upss) \).
Given \( \rups < \rrb \), let a monotonous function \( f(\rr,\rups) \) fulfill for some 
\( \rho < 1 \) 
\begin{EQA}[c]
	f(\rr,\rups) 
	\geq
	3 \nunu \rr \, \qqQ\bigl( \xx + \log(\rr/\rups) \bigr),
	\quad 
	\rups \leq \rr \leq \rrb,
\label{frhorrxxrr}
\end{EQA}
where the function \( \qqQ(\cdot) \) is given by \eqref{zzxxgfin}.
Then it holds
\begin{EQA}
    \P\biggl( 
        \sup_{\rups \leq \rr \leq \rrb} \,\,
        \sup_{\ups \in \Upss(\rr)} 
        \bigl\{ 
            \UP(\ups,\upss) - f\bigl( \rho^{-1} \rr, \rups \bigr)
        \bigr\} 
        \geq 0 
    \biggr)
    & \le &
    \frac{\rho}{1 - \rho} \ex^{-\xx} .
\label{Psuprzz}
\end{EQA}    
\end{theorem}

\begin{remark}
\label{RTsuprUP}
Formally the bound applies even with \( \rrb = \infty \) provided that \( (\CS{D}) \) 
is fulfilled on the whole set \( \Upsd \).
\end{remark}

\begin{remark}
If \( \gm = \infty \), then \( \qqQ(\xx) = \sqrt{2 \xx + 4 \dimp} \) and the condition \eqref{frhorrxxrr} on the drift simplifies to 
\( (3 \nunu \rr)^{-1} f(\rr,\rups) \geq \sqrt{2\xx + 4 \dimp + 2 \log(\rr/\rups)} \).
\end{remark}

\begin{proof}
By \eqref{frhorrxxrr} and Theorem~\ref{Tsmoothpenlc} for any \( \rr > \rups \)
\begin{EQA}
	&& \nquad
	\P\biggl(  
		\sup_{\ups \in \B_{\rr}(\upss) \setminus \B_{\rho \rr}(\upss)} 
        \bigl\{ \UP(\ups,\upss) - f\bigl( \rr,\rups) \bigr\}
		\geq 0 
	\biggr)
	\\
	& \leq &
	\P\biggl( 
		\frac{1}{3 \nunu \rr} \sup_{\ups \in \B_{\rr}(\upss)} \UP(\ups,\upss) 
		\geq 
		\qqQ\bigl(\xx + \log(\rr/\rups)\bigr) 
	\biggr)
	\leq 
	\frac{\rups}{\rr} \ex^{- \xx} .
\label{Puprrrho0}
\end{EQA}
Now define \( \rr_{k} = \rups \rho^{-k} \) for \( k=0,1,2,\ldots \).
Define also \( \kb \eqdef \log (\rrb/\rups) + 1 \).
It follows from \eqref{Puprrrho0} that 
\begin{EQA}
	&& \nquad
	\P\biggl( 
        \sup_{\ups \in \B_{\rrb}(\upss) \setminus \B_{\rups}(\upss)} 
        \biggl\{ 
            \UP(\ups,\upss) - f\bigl( \rho^{-1} \dist(\ups,\upss), \rups \bigr)
        \biggr\} \geq 0 
    \biggr)
    \\
    & \leq &
	\sum_{k=1}^{k^{*}} \P\biggl(  
		\frac{1}{\rr_{k}} \sup_{\ups \in \B_{\rr_{k}}(\upss) \setminus \B_{\rr_{k-1}}(\upss)} 
        \biggl\{ \UP(\ups,\upss) - f\bigl( \rr_{k},\rups \bigr) \biggr\}
		\geq 0 
	\biggr)
	\\
	& \leq &
	\ex^{-\xx} \sum_{k=1}^{\kb} \rho^{k}
	\leq 
	\frac{\rho}{1 - \rho} \ex^{- \xx} 
\label{PBrrbmBrups}
\end{EQA}
as required.
\end{proof}

\Section{A bound for the norm of a vector random process}
Let \( \UU(\ups) \), \( \ups \in \Ups \), be a smooth centered random vector process 
with values in \( \R^{\dimq} \), where \( \Ups \subseteq \R^{\dimp} \). 
Let also \( \UU(\upss) = 0 \) for a fixed point \( \upss \in \Ups \).
Without loss of generality assume \( \upss = 0 \).
We aim to bound the maximum of the norm \( \| \UU(\ups) \| \)
over a vicinity \( \Upss \) of \( \upss \).
Suppose that \( \UU(\ups) \) satisfies for each 
\( \gammav \in \R^{\dimp} \) and \( \alphav \in \R^{\dimq} \) with 
\( \| \gammav \| = \| \alphav \| = 1 \)
\begin{EQA}[c]
    \sup_{\ups \in \Ups} \log \E \exp\Bigl\{ 
		\lambda \gammav^{\T} \nabla \UU(\ups) \alphav 
	\Bigr\} 
    \le 
    \frac{\nunu^{2} \lambda^{2}}{2} , 
    \qquad 
    \lambda^{2} \leq 2 \gm^{2}.
\label{gUUgem}
\end{EQA}
Condition \eqref{gUUgem}
implies for any \( \ups \in \Upss \) with \( \| \ups \| \leq \rr \) and
\( \| \gammav \| = 1 \) in view of \( \UU(\upss) = 0 \)
\begin{EQA}[c]
    \log \E \exp\Bigl\{ 
    	\frac{\lambda}{\rr} \gammav^{\T} \UU(\ups) 
	\Bigr\} 
    \le 
    \frac{\nunu^{2} \lambda^{2} \| \ups \|^{2}}{2 \rr^{2}} ,
    \qquad 
    \lambda^{2} \leq 2 \gm^{2} ;
\label{gUUem}
\end{EQA}   
In what follows, we use the representation
\begin{EQA}[c]
    \| \UU(\ups) \|
    =
    \sup_{\| \uv \| \leq \rr} \frac{1}{\rr} \uv^{\T} \UU(\ups) .
\label{UU2ups}
\end{EQA}    
This implies for 
\( \Upss(\rr) = \bigl\{ \ups \in \Ups \colon \| \ups - \upss \| \leq \rr \bigr\} \)
\begin{EQA}[c]
    \sup_{\ups \in \Upss(\rr)} \| \UU(\ups) \|
    =
    \sup_{\ups \in \Upss(\rr)} \,\, \sup_{\| \uv \| \leq \rr} 
        \frac{1}{\rr}  \uv^{\T} \UU(\ups) .
\label{UU2vups}
\end{EQA}   
Consider a bivariate process \( \uv^{\T} \UU(\ups) \) of 
\( \uv \in \R^{\dimq} \) and \( \ups \in \Ups \subset \R^{\dimp} \).
By definition \( \E \uv^{\T} \UU(\ups) = 0 \).
Further, 
\( \nabla_{\uv} \bigl[ \uv^{\T} \UU(\ups) \bigr] = \UU(\ups) \) while 
\( \nabla_{\ups} \bigl[ \uv^{\T} \UU(\ups) \bigr] = \uv^{\T} \nabla \UU(\ups) 
= \| \uv \| \gammav^{\T} \nabla \UU(\ups) \) for \( \gammav = \uv / \| \uv \| \).
Suppose that \( \uv \in \R^{\dimq} \) and \( \ups \in \Ups \) are such that 
\( \| \uv \|^{2} + \| \ups \|^{2} \leq 2\rr^{2} \).
By the H\"older inequality, \eqref{gUUem}, and \eqref{gUUgem}, it holds for 
\( \| \gammav \| = \| \alphav \| = 1 \) and \( \ups \in \Upss(\rr) \)
\begin{EQA}
	&& \nquad
	\log \E \exp\biggl\{ 
		\frac{\lambda}{2 \rr} 
        (\gammav,\alphav)^{\T} \nabla \bigl[ \uv^{\T} \UU(\ups) \bigr] 
	\biggr\}
	\\
	& \leq &
	\frac{1}{2} \log \E \exp\biggl\{ 
		\frac{\lambda}{\rr} \gammav^{\T} \UU(\ups) 
	\biggr\}
	+
	\frac{1}{2} \log \E \exp\biggl\{ 
		\frac{\lambda}{\rr} \uv^{\T} \nabla \UU(\ups) \alphav
	\biggr\}
	\\
	& \leq &
	\frac{1}{2} \log \E \exp\biggl\{ 
		\frac{\lambda}{\rr} \gammav^{\T} \UU(\ups) 
	\biggr\} 
	+
	\frac{1}{2} \log \E \exp\biggl\{ 
		\frac{\lambda}{\rr} \| \uv \| \gammav^{\T} \nabla \UU(\ups) \alphav
	\biggr\} 
	\\
	& \leq &
	\frac{\nunu^{2} \lambda^{2}}{4 \rr^{2}} \bigl( \| \ups \|^{2} + \| \uv \|^{2} \bigr) 
	\leq
	\frac{\nunu^{2} \lambda^{2}}{2},
	\qquad
	|\lambda| \leq \gm .
\label{lexpuvUps}
\end{EQA}
We summarize our findings in the following theorem.

\begin{theorem}
\label{Tsqnorm}
Let a random \( \dimp \)-vector process \( \UU(\ups) \) for 
\( \ups \in \Ups \subseteq \R^{\dimp} \)
fulfill \( \UU(\upss) = 0 \), \( \E \UU(\ups) \equiv 0 \), 
and the condition \eqref{gUUgem} be satisfied. 
Then for each \( \rr \) and any \( \xx \ge 1/2 \), it holds 
\begin{EQA}
\label{Upsdboundno}
    \P \Bigl\{ 
        \sup_{\ups \in \Upss(\rr)} \bigl\| \UU(\ups) \bigr\| 
    	> 6 \nunu \rr \, \qqQ(\xx)
    \Bigr\}
    & \le & 
    \ex^{-\xx} ,
\end{EQA}
where \( \qqQ(\xx) \) is given by \eqref{zzxxgfin}.
\end{theorem}

\begin{proof}
The results follow from Theorem~\ref{Tsmoothpenlc} applied to the process 
\( \uv^{\T} \UU(\ups)/2 \) of the variable \( (\uv,\ups) \in \R^{\dimp + \dimq} \).
\end{proof}

\section{A deviation bound for the quadratic form \( \| \xiv \|^{2} \)}
This section presents a bound for a quadratic form \( \| \xiv \|^{2} \)
where \( \xiv = \DPc^{-1} \, \nabla \zeta(\thetavs) \).
The result only uses the condition \( (E\!D_{0}) \) which we restate in 
a slightly different form.
For \( \BB = \DPc^{-1} \VPc^{2} \DPc^{-1} \), define 
\begin{EQA}[c]
    \dimB
    \eqdef
    \tr \bigl( \BB \bigr) ,
    \qquad 
    \vA_{\BB}^{2}
    \eqdef
    2 \tr(\BB^{2}) ,
    \qquad
    \lambdaB \eqdef \lambda_{\max}\bigl( \BB \bigr). 
\label{BBrdd}
\end{EQA}   
Note that \( \dimB = \E \| \xiv \|^{2} \).
Moreover, if \( \xiv \) is a Gaussian vector then 
\( \vA_{\BB}^{2} = \Var\bigl( \| \xiv \|^{2} \bigr) \).
If \( \VPc^{2} = \DPc^{2} \), then \( \lambdaB = 1 \).
The condition \( (E\!D_{0}) \) means that 
the vector \( \VPc^{-1} \, \nabla \zeta(\thetavs) \) 
fulfills the following exponential moment condition:
\begin{EQA}[c]
    \log \E \exp\bigl( \gammav^{\T} \VPc^{-1} \, \nabla \zeta(\thetavs) \bigr) 
    \le 
    \| \gammav \|^{2}/2,
    \qquad 
    \gammav \in \R^{\dimp}, \, \| \gammav \| \le \gm .
\label{expgamgm}
\end{EQA}
Here \( \nunu \) is set to one.
\cite{SP2011} argued how the case of any \( \nunu \ge 1 \) can be reduced to 
\( \nunu \approx 1 \) by a slight change of scale and reducing the value \( \gm \) which 
is typically large.
For ease of presentation, suppose that \( \gm^{2} \ge 2 \dimB \).
The other case only changes the constants in the inequalities. 
Note that \( \| \xiv \|^{2} = \etav^{\T} \BB \, \etav \).
Define \( \muc = 2/3 \) and
\begin{EQA}
    \gmc
    & \eqdef &
    \sqrt{\gm^{2} - \muc \dimB} ,
    \\
    2\xxc
    & \eqdef &
    (\gm^{2}/\muc - \dimB)/\lambdaB + \log \det \bigl( \Id_{\dimp} - \muc \BB/\lambdaB \bigr) .
\label{yycgmcxxcAd}
\end{EQA}

\begin{proposition}[\cite{SP2011}]
\label{LLbrevelocro}   
Let \( (E\!D_{0}) \) hold with \( \nunu = 1 \) and 
\( \gmb^{2} \ge 2 \dimB \).
Then
for each \( \xx > 0 \)
\begin{EQA}
    \P\bigl( \| \xiv \| \ge \qq(\BB,\xx) \bigr)
    & \le &
    2 \ex^{-\xx} + 8.4 \ex^{-\xxc} ,
\label{PxivbzzBBro}
\end{EQA}    
where \( \qq(\BB,\xx) \) is defined by
\begin{EQA}
\label{PzzxxpBro}
    \qq^{2}(\BB,\xx)
    & \eqdef &
    \begin{cases}
        \dimB + 2 \vA_{\BB} \xx^{1/2}, &  \xx \le \vA_{\BB}/(18 \lambdaB) , \\
        \dimB + 6 \lambdaB \xx, & \vA_{\BB}/(18 \lambdaB) < \xx \le \xxc , \\
        \bigl| \yyc + 2 \lambdaB (\xx - \xxc)/\gmc \bigr|^{2}, & \xx > \xxc .
    \end{cases}
\label{zzxxppdBlro}
\end{EQA}    
with \( \yyc^{2} \leq \dimB + 6 \lambdaB \xxc \).
\end{proposition}

Depending on the value \( \xx \), we observe three types of tail behavior of the 
quadratic form \( \| \xiv \|^{2} \). 
The sub-Gaussian regime for 
\( \xx \le \vA_{\BB}/(18 \lambdaB) \) 
and the Poissonian regime for \( \xx \le \xxc \) 
are similar to the case of a Gaussian quadratic form.
The value \( \xxc \) from \eqref{yycgmcxxcAd} is of order \( \gm^{2} \).
In all our results we suppose that \( \gm^{2} \) and hence, \( \xxc \) 
is sufficiently large and 
the quadratic form \( \| \xiv \|^{2} \) can be bounded with a dominating 
probability by \( \dimB + 6 \lambdaB \xx \) for a proper \( \xx \).
%
We refer to \cite{SP2011} for the proof of this and related results, further 
discussion and references.

\section{Some inequalities for the normal law}
This section collects some simple but useful facts about the properties of the multivariate standard 
normal distribution.
Many similar results can be found in the literature, we present the proofs to keep the presentation 
self-contained. 
Everywhere in this section \( \gammav \) means a standard normal vector in \( \R^{\dimp} \).

\begin{lemma}
\label{LGaussmom}
Let \( \mu \in (0,1) \).
Then for any vector \( \lambdav \in \R^{\dimp} \) with
\( \| \lambdav \|^{2} \le \dimp \) and any \( \rr > 0 \)
\begin{EQA}[c]
    \log \E \bigl\{
        \exp (\lambdav^{\T} \gammav) \Ind\bigl( \| \gammav \| > \rr \bigr)
    \bigr\}
    \le
    - \frac{1 - \mu}{2} \rr^{2} + \frac{1}{2 \mu} \| \lambdav \|^{2}
    + \frac{\dimp}{2} \log(\mu^{-1}) .
    \qquad
\label{lamgamrr}
\end{EQA}
Moreover, if \( \rr^{2} \ge 6 \dimp + 4 \xx \), then 
\begin{EQA}[c]
    \E \Bigl\{
        \exp (\lambdav^{\T} \gammav) \Ind\bigl( \| \gammav \| \le \rr \bigr)
    \Bigr\}
    \ge
    \ex^{ \| \lambdav \|^{2}/2} \bigl( 1 - \ex^{-\xx} \bigr) .
\label{exIdgamrr}
\end{EQA}
\end{lemma}

\begin{proof}
We use that for \( \mu < 1 \)
\begin{EQA}[c]
    \E \bigl\{
        \exp(\lambdav^{\T} \gammav) \Ind\bigl( \| \gammav \| > \rr \bigr)
    \bigr\}
    \le
    \ex^{- (1 - \mu) \rr^{2}/2}
        \E \exp \bigl\{ \lambdav^{\T} \gammav + (1 - \mu) \| \gammav \|^{2}/2 \bigr\} .
\label{explamgamrr}
\end{EQA}
It holds
\begin{EQA}
    \E \exp \bigl\{ \lambdav^{\T} \gammav + (1 - \mu) \| \gammav \|^{2}/2 \bigr\}
    &=&
    (2 \pi)^{-\dimp/2} \int
        \exp \bigl\{ \lambdav^{\T} \gammav - \mu \| \gammav \|^{2}/2 \bigr\} d \gammav
    \\
    &=&
    \mu^{-\dimp/2} \exp\bigl( \mu^{-1} \| \lambdav \|^{2}/2 \bigr)
\label{exp1mummu}
\end{EQA}
and \eqref{lamgamrr} follows.

Now we apply this result with \( \mu = 1/2 \).
In view of
\( \E \exp(\lambdav^{\T} \gammav) = \ex^{\| \lambdav \|^{2}/2} \),
\( \rr^{2} \geq 6 \dimp + 4 \xx \),
and \( 2 + \log(2) < 3 \), it follows
for \( \| \lambdav \|^{2} \le \dimp \)
\begin{EQA}
    && \nquad
    \ex^{- \| \lambdav \|^{2}/2}
    \E \bigl\{
        \exp (\lambdav^{\T} \gammav) \Ind\bigl( \| \gammav \| \le \rr \bigr)
    \bigr\}
    \\
    & \ge &
    1 - \exp\bigl( - \rr^{2}/4 + \dimp + (\dimp/2) \log(2) \bigr)
    \ge
    1 - \exp(- \xx)
\label{exIdgamrrd}
\end{EQA}
which implies \eqref{exIdgamrr}.
\end{proof}

\begin{lemma}
For any \( \uv \in \R^{\dimp} \),
any unit vector \( \av \in \R^{\dimp} \), and any \( \qq > 0 \), it holds
\begin{EQA}
\label{gampxiz12}
	\P\bigl( \| \gammav - \uv \| \geq \qq \bigr)
	& \leq &
	\exp\bigl\{ - \qq^{2}/4 + \dimp/2 + \| \uv \|^{2}/2 \bigr\} ,
	\\
	\E \bigl\{ | \gammav^{\T} \av |^{2} \Ind\bigl( \| \gammav - \uv \| \geq \qq \bigr) \bigr\}
	& \leq &
	(2 + |\uv^{\T} \av|^{2}) \exp\bigl\{ - \qq^{2}/4 + \dimp/2 + \| \uv \|^{2}/2 \bigr\} .
\label{gampxiz122}
\end{EQA}
\end{lemma}

\begin{proof}
By the exponential Chebyshev inequality, for any \( \lambda < 1 \) 
\begin{EQA}
	\P\bigl( \| \gammav - \uv \| \geq \qq \bigr)
	& \leq &
	\exp \bigl( - \lambda \qq^{2}/2 \bigr) \E \exp\bigl( \lambda \| \gammav - \uv \|^{2}/2 \bigr)
	\\
	&=&
	\exp\Bigl\{ 
		- \frac{\lambda \qq^{2}}{2} - \frac{\dimp}{2} \log(1 - \lambda) 
		+ \frac{\lambda}{2 (1-\lambda)} \| \uv \|^{2}
	\Bigr\} .
\label{gampxiz}
\end{EQA}
In particular, with \( \lambda = 1/2 \), this implies \eqref{gampxiz12}.
Further, for \( \| \av \| = 1 \)
\begin{EQA}
	\E \bigl\{ | \gammav^{\T} \av |^{2} \Ind(\| \gammav - \uv \| \geq \qq) \bigr\}
	& \leq &
	\exp \bigl( - \qq^{2}/4 \bigr) 
	\E \bigl\{ | \gammav^{\T} \av |^{2} \exp\bigl( \| \gammav - \uv \|^{2}/4 \bigr) \bigr\}
	\\
	& \leq &
	(2 + | \uv^{\T} \av |^{2}) \exp \bigl( - \qq^{2}/4 + \dimp/2 + \| \uv \|^{2}/2 \bigr) 
\label{•}
\end{EQA}
and \eqref{gampxiz122} follows.
\end{proof}

The next result explains the concentration effect for the norm \( \| \xiv \|^{2} \)
of a Gaussian vector.
We use a version from \cite{SP2011}.

\begin{lemma}
\label{LPGxx}
For each \( \xx \),
\begin{EQA}[rclcl]
\label{Pgamlxx}
	\P\bigl( \| \gammav \| \geq \qq(\dimp,\xx) \bigr)
	& \leq &
	\exp\bigl( - \xx \bigr) ,
	\quad
	\P\bigl( \| \gammav \| \leq \qq_{1}(\dimp,\xx) \bigr)
	& \leq &
	\exp\bigl( - \xx \bigr) ,
	\qquad
\label{Pgammxx}
\end{EQA}
where 
\begin{EQA}
	\qq^{2}(\dimp,\xx) 
	& \eqdef &
	\dimp + \sqrt{6.6 \dimp \xx} \vee (6.6 \xx) ,
	\qquad
	\qq_{1}^{2}(\dimp,\xx) \eqdef \dimp - 2 \sqrt{\dimp \xx} .
\label{zz0xxG}
\end{EQA}
\end{lemma}

The next lemma bounds from above the Kullback-Leibler divergence between 
two normal distributions.

\begin{lemma}
\label{KullbTVd}
Let \( \P_{0} = \ND(0,\Id_{\dimp}) \) and \( \P_{1} = \ND(\betav,(\DD^{\T} \DD)^{-1}) \)
some non-degenerated matrix \( \DD \).
If 
\begin{EQA}
	\| \DD^{\T} \DD - \Id_{\dimp} \|_{\infty} 
	& \leq &
	\rd \leq 1/2 ,
\label{DPcbIdeps}
\end{EQA}
then 
\begin{EQA}
	2 \kullb(\P_{0},\P_{1})
	&=&
	- 2 \E_{0} \log \frac{d\P_{1}}{d\P_{0}}
	\\
	& \leq &
 	\tr (\DD^{\T} \DD - \Id_{\dimp})^{2} + (1 + \rd) \| \betav \|^{2}
	\leq  
	\rd^{2} \, \dimp + (1 + \rd) \| \betav \|^{2} .
\label{kullbP0P1}
\end{EQA}
For any measurable set \( A \subset \R^{\dimp} \), it holds with 
\( \gammav \sim \ND(0,\Id_{\dimp}) \)
\begin{EQA}[c]
	\bigl| \P_{0}(A) - \P_{1}(A) \bigr|
	=
	\bigl| \P\bigl( \gammav \in A \bigr)
	- \P\bigl( \DD (\gammav - \betav) \in A \bigr) \bigr| 
	\leq 
	\sqrt{\kullb(\P_{0},\P_{1}) / 2}.
\end{EQA}
\end{lemma}

\begin{proof}
It holds 
\begin{EQA}[c]
	2 \log \frac{d\P_{1}}{d\P_{0}}(\gammav)
	=
	\log \det (\DD^{\T} \DD) 
	- (\gammav - \betav)^{\T} \DD^{\T} \DD (\gammav - \betav)  
	+ \| \gammav \|^{2}
\end{EQA}
with \( \gammav \) standard normal and
\begin{EQA}
	2 \kullb(\P_{0},\P_{1})
	&=&
	- 2 \E_{0} \log \frac{d\P_{1}}{d\P_{0}}
	=
	- \log \det (\DD^{\T} \DD) 
	+ \tr (\DD^{\T} \DD - \Id_{\dimp}) 
	+ \betav^{\T} \DD^{\T} \DD \betav .
\end{EQA}
Let \( a_{j} \) be the \( j \)th eigenvalue of \( \DD^{\T} \DD - \Id_{\dimp} \).  
\( \| \DD^{\T} \DD - \Id_{\dimp} \|_{\infty} \leq \rd \leq 1/2 \) yields 
\( |a_{j}| \le 1/2 \) and 
\begin{EQA}
	2 \kullb(\P_{0},\P_{1})
	&=&
    \betav^{\T} \DD^{\T} \DD \betav 
	+
    \sum_{j=1}^{\dimp} \bigl\{ a_{j} - \log(1 + a_{j}) \bigr\}
	\\
	& \leq &
    (1 + \rd) \| \betav \|^{2} 
	+ \sum_{j=1}^{\dimp} a_{j}^{2} 
	\\
	& \leq &
	(1 + \rd) \| \betav \|^{2} + \tr (\DD^{\T} \DD - \Id_{\dimp})^{2} 
	\leq 
	(1 + \rd) \| \betav \|^{2} + \rd^{2} \, \dimp.
\end{EQA}
This implies by Pinsker's inequality 
\begin{EQA}
	\sup_{A} | \P_{0}(A) - \P_{1}(A) |
	& \leq &
	\sqrt{\frac{1}{2} \kullb(\P_{0},\P_{1})} 
\label{TVdPiBvM}
\end{EQA}
as required.
\end{proof}

\bibliography{exp_ts,listpubm-with-url}

\end{document}